\documentclass[11pt,reqno]{amsart}
\usepackage[a4paper,margin=28mm]{geometry}
\usepackage{amsmath,amssymb,mathtools}
\usepackage{tikz}
\usepackage[colorlinks=true,linkcolor=blue,citecolor=blue,urlcolor=blue]{hyperref}
\numberwithin{equation}{section}
\newtheorem{theorem}{Theorem}[section]
\newtheorem{lemma}[theorem]{Lemma}
\newtheorem{proposition}[theorem]{Proposition}
\newtheorem{corollary}[theorem]{Corollary}
\theoremstyle{remark}

\newcommand{\Pp}{\mathbb P}
\newcommand{\E}{\mathbb E}
\newcommand{\dist}{\operatorname{dist}}
\newcommand{\diam}{\operatorname{diam}}
\newcommand{\Pois}{\operatorname{Poisson}}

\title[Dimension-uniform spherical connectivity]
{Dimension-uniform connectivity laws and hitting times
for spherical random geometric graphs}
\author{Wentao Hu}
\address{Qiuzhen College, Tsinghua University, Beijing, China}
\email{huwt24@mails.tsinghua.edu.cn}
\begin{document}
\begin{abstract}
Let $G_d(n,p)$ be the graph formed by joining pairs of $n$ independent
uniform points on $\mathbb S^{d-1}$ within a radius chosen to give edge
probability $p$. For each fixed integer $K\ge0$ and $t\in\mathbb R$,
we prove that at $np=\log n+K\log\log n+t$, the probabilities of
$(K+1)$-vertex-connectivity and $(K+1)$-edge-connectivity both converge
 to $\exp(-e^{-t}/K!)$ uniformly over all $d=d(n)\ge3$. This is the same limit as in
the Erd\H{o}s--R\'enyi graph. Moreover, as the radius increases on a single sample, the first radii for
minimum degree at least $K+1$, $(K+1)$-vertex-connectivity, and
$(K+1)$-edge-connectivity coincide with probability tending to one
uniformly over $d\ge3$.
To our knowledge, this is the first critical-window law for connectivity
in a random geometric graph that holds uniformly from fixed dimensions
to dimensions growing arbitrarily with $n$.
\end{abstract}

\maketitle

\section{Introduction}\label{sec:introduction}

Connectivity in the Erd\H{o}s--R\'enyi random graph $G(n,p)$ is governed by
isolated vertices. At $p=(\log n+t)/n$, with $t$ fixed, the number of
isolated vertices converges to a Poisson variable of mean
$e^{-t}$, and with high probability all other vertices lie in one
component. The limiting probability of connectivity is therefore
$\exp(-e^{-t})$. For $k$-connectivity the corresponding obstruction is
a vertex of degree less than $k$. The classical theorem of
Erd\H{o}s and R\'enyi~\cite{ErdosRenyi1961}, originally proved for a
prescribed number of edges, places the transition for fixed $k$ in the
window $np=\log n+(k-1)\log\log n+O(1)$. Bollob\'as and Thomason~\cite{BollobasThomason1985} proved the
corresponding hitting-time result for the random graph process. When the edges are added in uniformly
random order, with high probability the graph first becomes
$k$-vertex-connected at the step when its minimum degree reaches $k$.
The first step giving $k$-edge-connectivity lies between these two and
hence coincides with them.

Similar connectivity phenomena arise in random geometric graphs of
fixed dimension, even though nearby vertices tend to form a small
component in which every vertex has the required degree. For uniform points in a cube of fixed dimension at least two,
Penrose~\cite{Penrose1999} showed that, with high probability, the first
radius giving minimum degree $k$ also gives $k$-vertex-connectivity.
Penrose and Yang~\cite{PenroseYang2025} determined the fluctuations of
these radii in Euclidean regions, where points near the boundary can
govern the threshold. On a torus of fixed dimension at least two, Ignasiak and
Lichev~\cite{IgnasiakLichev2025} studied a related process in which
points are added one at a time at a fixed connection radius $r$.
They showed that, with probability tending to one as $r\to0$, the
graph remains $k$-connected from the first time its minimum degree
reaches $k$ onward.

This correspondence provides a way to determine limiting
connectivity probabilities through the distribution of
low-degree vertices.
For example,  on a flat torus of fixed dimension at least two,
Penrose's results~\cite{Penrose2003} show that, for each fixed $k\ge1$,
the critical window for $k$-vertex-connectivity in the edge probability
$p$ and the limiting probability throughout that window agree with
those of $G(n,p)$.

The picture is less complete when the dimension grows with $n$.
In sufficiently high dimensions, Bangachev and
Bresler~\cite{BangachevBresler2024} construct a coupling that,
with high probability, sandwiches the spherical geometric graph
between two Erd\H{o}s--R\'enyi graphs with slightly smaller and
larger edge probabilities.
At the connectivity scale $np\sim\log n$, the corresponding
changes in $np$ are $O((\log n)^{7/2}/\sqrt d)$.
Thus, when $d/(\log n)^7\to\infty$, these changes are $o(1)$,
and the coupling transfers the Erd\H{o}s--R\'enyi limiting
$k$-connectivity probability throughout the critical window.
However, this comparison leaves the critical-window behaviour
unresolved when the dimension grows slowly with $n$.
To our knowledge, no previous work has determined a
connectivity critical-window law for a standard random
geometric graph model uniformly from fixed dimensions through
arbitrary dimension growth.

We resolve this question for the spherical geometric graph $G_d(n,p)$,
formed by joining points within distance $r$ in a sample of $n$
independent uniform points on $\mathbb S^{d-1}$. Our results hold for
every fixed $k\ge1$, uniformly over all dimensions $d=d(n)\ge3$.
In terms of the edge probability $p$, the critical window and the
limiting probabilities of $k$-vertex-connectivity
throughout that window agree with those of $G(n,p)$. The corresponding
hitting-time law holds as well. As the radius increases on a single
sample, the first radius at which the minimum degree reaches $k$ is,
with high probability, also the first radius for both
$k$-vertex-connectivity and $k$-edge-connectivity.

\subsection{Model and main results}

Let $X_1,\ldots,X_n$ be independent uniform points on
$\mathbb S^{d-1}$, where $d\ge3$, and write
$\dist(x,y)=\|x-y\|$ for chordal distance. For $0<p<1$, let
$r=r_d(p)$ be the unique radius satisfying
$\Pp\{\dist(X_1,X_2)\le r\}=p$. The graph $G_d(n,p)$ has vertex
set $[n]$ and joins $i$ and $j$ whenever $\dist(X_i,X_j)\le r$.

Fix an integer \(K\ge0\). Write $\delta(G)$ for
the minimum degree of $G$. A graph on more than $K+1$ vertices is
$(K+1)$-vertex-connected if it remains connected after any deletion
of at most $K$ vertices, and $(K+1)$-edge-connected if it remains
connected after any deletion of at most $K$ edges.

By rotational invariance, conditional on $X_i$, each of the other
$n-1$ points is independently adjacent to $X_i$ with probability
$p$. Thus the expected number of vertices of degree $K$ is
\begin{equation}\label{eq:Lambda}
 \lambda_{K,n}
 =n\binom{n-1}{K}p^K(1-p)^{n-1-K}.
\end{equation}

\begin{theorem}[Connectivity and minimum degree]\label{thm:main}
Fix $0\le T<\infty$ and suppose that
\[
 np=\log n+K\log\log n+t,\qquad |t|\le T.
\]
Let $\Phi_K$  denote either the event that $G=G_d(n,p)$ is
$(K+1)$-vertex-connected or the event that it is $(K+1)$-edge-connected.
For all sufficiently large $n$, uniformly over $d\ge3$ and $|t|\le T$,
\begin{align}
 \Pp\bigl(\{\delta(G)\ge K+1\}\cap\Phi_K^c\bigr)
 &\le \frac{C_{K,T}}{\log n},\label{eq:main-equivalence}\\
 \bigl|\Pp(\Phi_K)-e^{-\lambda_{K,n}}\bigr|
 &\le \frac{C_{K,T}}{\log n}.\label{eq:main-probability}
\end{align}
In particular, if $np=\log n+K\log\log n+t+o(1)$ for fixed
$t\in\mathbb R$, then for either event, uniformly over $d\ge3$,
\[
 \Pp(\Phi_K)\longrightarrow
 \exp\!\left(-\frac{e^{-t}}{K!}\right).
\]
\end{theorem}

For $K=0$, the proof also describes the graph in the connectivity
window. Uniformly over $d\ge3$ and $|np-\log n|\le T$,
\begin{equation}\label{eq:ordinary-structure}
 \Pp\{G_d(n,p)\text{ has exactly one component of order at least two}\}
 \ge 1-\frac{C_T}{\log n}.
\end{equation}

Now keep the sample $X_1,\ldots,X_n$ fixed and increase $p$.
Let $p_{\delta,K}$, $p_{\mathrm e,K}$, and $p_{\kappa,K}$
be the first values of $p$ at which the graph has minimum degree
at least $K+1$, is $(K+1)$-edge-connected, and is
$(K+1)$-vertex-connected, respectively. Since $r_d(p)$ is strictly
increasing, these parameters coincide exactly when the corresponding
radii coincide.

\begin{theorem}[Coincident hitting times]\label{thm:hitting}
For every fixed integer $K\ge0$,
\[
 \inf_{d\ge3}
 \Pp\{p_{\delta,K}=p_{\mathrm e,K}=p_{\kappa,K}\}
 \longrightarrow 1.
\]
\end{theorem}

\begin{corollary}[Zero--one law for connectivity]\label{cor:zero-one}
For any sequence $p_n\in(0,1)$, uniformly over $d\ge3$,
\[
 \Pp\{G_d(n,p_n)\text{ is connected}\}
 \longrightarrow
 \begin{cases}
  0, & np_n-\log n\to-\infty,\\
  1, & np_n-\log n\to+\infty.
 \end{cases}
\]
\end{corollary}

\subsection{Related work}\label{sec:related-work}

Our results cover all dimensions $d\ge3$. In dimension $d=2$,
the underlying space is the circle $\mathbb S^1$, where the two
thresholds differ. Maehara~\cite{Maehara1990} showed that, to first
order, isolated vertices disappear at $p=(\log n)/n$, while
connectivity requires twice that edge probability.

On $\mathbb S^2$, Maehara~\cite{Maehara2004} showed that when
$np=\log n+t+o(1)$ for fixed $t$, the number of isolated vertices
in $G_3(n,p)$ converges to a Poisson variable of mean $e^{-t}$.
For $p=c(\log n)/n$, Maehara asked for the connectivity threshold
constant $c$, which his bounds placed between $1$ and $2$
\cite[Problem~1]{Maehara2004}. The connectivity bound of Song,
Goeckel, and Towsley~\cite[Theorem~2]{SongGoeckelTowsley2006},
together with Maehara's lower bound, gives $c=1$.
Theorem~\ref{thm:main} also gives the limiting connectivity
probability throughout the critical window, uniformly over $d\ge3$

For each fixed $d\ge3$, the spherical annulus estimates of Galhotra,
Mazumdar, Pal, and Saha~\cite{GalhotraMazumdarPalSaha2019} imply
that $G_d(n,p)$ is connected with high probability whenever
$p\ge C_d(\log n)/n$ for a sufficiently large constant $C_d$.
David and Feige~\cite{DavidFeige2013} fixed an angular connection radius $\theta\in(0,\pi/2)$ and let $d$ grow. Writing $p$ for the resulting edge probability,
they showed that the graph is connected with high probability once
$n\ge C p^{-1}\log(1/p)$, a constant multiple of the sample size
needed for isolated vertices to disappear.

Narayanan~\cite[Theorem~6.4]{Narayanan2026} proved that increasing
graph properties have thresholds in the connection radius on spheres,
tori, and cubes even when the dimension varies with $n$. He asked
whether such thresholds also exist when the parameter is the edge
probability~\cite[Question~7.4]{Narayanan2026}.
Theorem~\ref{thm:main} answers this question for vertex- and
edge-connectivity of each fixed order on spheres with $d\ge3$,
and gives the limiting probability throughout the critical window.

\subsection{Proof strategy and organization}

In the critical window, we show that with probability
$1-O_{K,T}(1/\log n)$, every deletion of at most $K$ vertices leaves
exactly one component of order at least two. If $\delta(G)\ge K+1$,
every surviving vertex retains a neighbor and hence belongs to that
component. This gives $(K+1)$-vertex-connectivity, and therefore also
$(K+1)$-edge-connectivity.

In Sections~\ref{sec:geometry} and~\ref{sec:local}, we show that,
with high probability, every nontrivial component remaining
after any deletion of at most $K$ vertices has more than
$\lceil(\log n)^8\rceil$ vertices.

In Section~\ref{sec:structure}, we split the sample into a base set
and a reserve set. With high probability, after any deletion of at
most $K$ vertices, each component of the base graph is either an
isolated vertex or a large component. Adding the surviving reserve
points one by one then joins all the large components into one. Too few vertices remain
outside it to form another large component, so the preceding estimate
shows that they are all isolated.

In Section~\ref{sec:consequences}, we use Poissonization to study the
distribution of the number of vertices of degree exactly $K$.
The resulting estimate for the probability that this count is zero,
together with the structural result, gives the limiting connectivity
probability. We then show that the three hitting times coincide
with high probability

\section{Dimension-uniform cap geometry}\label{sec:geometry}

Let $\mu$ be normalized surface measure on $\mathbb S^{d-1}$. For
$x\in\mathbb S^{d-1}$ and $0\le s\le2$, let
\[
 B(x,s)=\{y\in\mathbb S^{d-1}:\dist(x,y)\le s\}
\]
be the closed chordal ball of radius $s$ centered at $x$.
By rotational invariance, its mass is independent of $x$.
Write $F(s)=\mu(B(x,s))$, so that the connection radius $r=r_d(p)$
satisfies $F(r)=p$.

Put $\tau=1-r^2/2$. Since $\dist(x,y)^2=2(1-x\cdot y)$, the
neighborhood of $x$ is the spherical cap
\[
 C_x=B(x,r)=\{y\in\mathbb S^{d-1}:x\cdot y\ge\tau\}.
\]
Write $q(x,y)=\mu(C_x\cap C_y)$ for the mass of the overlap.

In the critical window of Theorem~\ref{thm:main}, all asymptotic
estimates are uniform over the indicated dimensions and $|t|\le T$,
with constants that may depend on $K$ and $T$.

Now fix $x\in\mathbb S^{d-1}$ and let $Y$ have distribution $\mu$.
By rotational invariance, $x\cdot Y$ has the same distribution as
the first coordinate $Y_1$. 
Thus, for $-1\le u\le 1$, 
\[
\begin{aligned}
 \Pp\{x\cdot Y\le u\}
 &=\Pp\{Y_1\le u\}\\
 &=\frac{1}{|\mathbb S^{d-1}|}
   \int_{-1}^{u}|\mathbb S^{d-2}|\,
   (1-z^2)^{(d-2)/2}\frac{dz}{\sqrt{1-z^2}}\\
 &=a_d\int_{-1}^{u}(1-z^2)^{(d-3)/2}\,dz,
\end{aligned}
\]
where
\[
 a_d:=\frac{|\mathbb S^{d-2}|}{|\mathbb S^{d-1}|}
 =\frac{\Gamma(d/2)}{\sqrt{\pi}\,\Gamma((d-1)/2)}
 \asymp\sqrt d.
\]
Write
\[
 J(u)=\int_u^1(1-z^2)^{(d-3)/2}\,dz,\qquad -1\le u\le1.
\]
One has
\[
 F(s)=\Pp\{x\cdot Y\ge1-s^2/2\}=a_dJ(1-s^2/2),
 \qquad 0\le s\le2.
\]
In particular, 
\[
p=F(r)=a_dJ(\tau).
\]
By symmetry, $F(\sqrt2)=1/2$, so $p\le1/2$ implies $r\le\sqrt2$
and $\tau\ge0$. In the window above, $p<1/2$ for sufficiently large
$n$, and hence $r<\sqrt2$ and $\tau>0$.

Differentiating the tail formula for $F(s)$ gives
\begin{equation}\label{eq:chord-density}
 F'(s)=a_d s^{d-2}(1-s^2/4)^{(d-3)/2},\qquad 0<s<2.
\end{equation}

Moreover, let $B^{d-1}(0,s)$ denote the Euclidean ball of radius $s$ in
$\mathbb R^{d-1}$ and set
\[
 w_d:=\frac{a_d}{d-1}.
\]
Since $|B^{d-1}(0,1)|=|\mathbb S^{d-2}|/(d-1)$, one has
\begin{equation}\label{eq:geometry_for_w_d}
     \frac{|B^{d-1}(0,s)|}{|\mathbb S^{d-1}|}=w_d s^{d-1}.
\end{equation}

\begin{lemma}[Mean cap overlap]\label{lem:overlap}
For every $d\ge2$ and $x\in\mathbb S^{d-1}$,
\[
 \int q(x,y)\,\mu(dy)=p^2.
\]
\end{lemma}
\begin{proof}
Fix $x$ and sample $Y,Z$ independently from $\mu$.
Conditional on $Y=y$, the probability that $Z\in C_x\cap C_Y$
is $q(x,y)$. Conditional on $Z=z$, the event $z\in C_Y$
is equivalent to $Y\in C_z$ and has probability $p$ by independence.
Thus
\[
 \int q(x,y)\,\mu(dy)
 =\Pp\{Z\in C_x,\ Y\in C_Z\}
 =p^2. \qedhere
\]
\end{proof}

\begin{lemma}[Radius and mass of a single cap]\label{lem:cap-estimates}
Suppose $d\ge3$ and $(\log n)/(2n)\le p\le2\log n/n$.
For all sufficiently large $n$, uniformly in $d$ and $p$,
\begin{align}
 r&\asymp(\log n/n)^{1/(d-1)},\label{eq:radius-scale}\\
 p_\lambda:=F(\lambda r)
 &\le 2^{(d-3)/2}\lambda^{d-1}p
 \le(\sqrt2\,\lambda)^{d-1}p,\qquad 0\le\lambda r\le2.
 \label{eq:near-cap}
\end{align}
\end{lemma}

\begin{proof}
Integrating the density \eqref{eq:chord-density} gives
\[
 p=F(r)=a_d\int_0^r s^{d-2}(1-s^2/4)^{(d-3)/2}\,ds.
\]
For $0\le s\le r\le\sqrt2$, the factor $(1-s^2/4)^{(d-3)/2}$
lies between $2^{-(d-3)/2}$ and $1$. Since
$a_d\int_0^r s^{d-2}\,ds=w_d r^{d-1}$, we obtain
\[
 2^{-(d-3)/2}w_d r^{d-1}\le p\le w_d r^{d-1}.
\]
Using $w_d\asymp d^{-1/2}$ and taking $(d-1)$st roots proves
\eqref{eq:radius-scale}.

For $\lambda\ge0$ with $\lambda r\le2$, the same density bound gives
\[
 F(\lambda r)
 \le w_d(\lambda r)^{d-1}
 \le 2^{(d-3)/2}\lambda^{d-1}p
 \le(\sqrt2\,\lambda)^{d-1}p.
\]
This proves \eqref{eq:near-cap}. In particular, taking $\lambda=1/10$
gives $F(r/10)\le Ce^{-d}p$.
\end{proof}

\begin{lemma}[Union and intersection of two caps]\label{lem:crescent}
For $d\ge2$ and $0<p\le1/2$, put $u=\dist(x,y)/r$. Then
\begin{equation}\label{eq:cap-union}
\mu(C_x\cup C_y)\ge p\left(1+\frac{u\wedge 1}{\pi}\right).
\end{equation}
Moreover, there are absolute constants $a,C>0$ such that, for all sufficiently
large $n$, all $d\ge3$, and all $(\log n)/(2n)\le p\le 2\log n/n$,
\begin{equation}\label{eq:far-overlap}
 \dist(x,y)\ge r/10
 \quad\Longrightarrow\quad
 q(x,y)\le C e^{-a(d\wedge\log n)}p.
\end{equation}
\end{lemma}

\begin{lemma}[Spherical neighborhood growth]
\label{lem:spherical-growth}
There is an absolute constant $c>0$ such that the following
holds for every $d\ge3$. Let $B\subseteq\mathbb S^{d-1}$ be
measurable with $0<\mu(B)\le1/2$, and define its closed chordal
$h$-neighborhood by
\begin{equation}\label{def:B_h}
 B_h=\left\{z\in\mathbb S^{d-1}:
          \inf_{w\in B}\dist(z,w)\le h\right\}.
\end{equation}
If $0<h\sqrt{d-1}\le 1$, then
\begin{equation}\label{eq:spherical-growth}
 \mu(B_h)-\mu(B)
 \ge c h\sqrt{d-1}\,\mu(B)^{1-1/(d-1)}.
\end{equation}
\end{lemma}

The proofs of Lemmas~\ref{lem:crescent} and~\ref{lem:spherical-growth}
are given in Appendices~\ref{app:two-cap-proof} and~\ref{app:growth-proof},
respectively.

\section{Uniform exclusion of small separated sets}\label{sec:local}

Throughout this section set $b=\lceil(\log n)^8\rceil$, $\delta=1/10$.
A vertex set $S$ in a graph $H$ is \emph{$K$-separated} if
$H[S]$ is connected and its external vertex neighborhood
\[
 \partial_H S=\{v\notin S:\text{$v$ has a neighbor in $S$}\}
\]
has order at most $K$. In particular, every nontrivial component left
after deleting at most $K$ vertices is a $K$-separated set in the
original graph.

We show that small nontrivial $K$-separated sets are absent with high
probability, uniformly over $d\ge3$.

\begin{proposition}\label{prop:local}
There is a constant $C_{K,T}$ such that, for all sufficiently large $n$,
uniformly over $d\ge3$ and $|t|\le T$, the probability that $G_d(n,p)$
contains a $K$-separated set of order $2,\ldots,b$ is at most
$C_{K,T}/\log n$.
\end{proposition}

We use repeatedly the uniform baseline estimates
\begin{equation}\label{eq:local-baseline}
 np\asymp\log n,\qquad
 (\log n)^K n e^{-np}=O_{K,T}(1)
\end{equation}
and the elementary binomial bound
\begin{equation}\label{eq:local-binomial}
 \Pp\{\operatorname{Bin}(N,v)\le K\}
 \le C_K(1+Nv)^K e^{-Nv}.
\end{equation}

\subsection{Sets of very small diameter}

Let $\mathcal A$ be the event that $G_d(n,p)$
contains a $K$-separated set $S$ with $|S|\ge2$ and
$\diam(S)\le\delta r$.

For distinct points $x,y$ with $s=\dist(x,y)\le\delta r$, define
\[
 R(x,y)=(C_x\cup C_y)\setminus B(x,s).
\]
For each ordered pair of distinct labels $i,j$, let
\[
 E_{ij}=
 \left\{
 \begin{aligned}
 &0<\dist(X_i,X_j)\le\delta r,\\
 &\#\{k\notin\{i,j\}:X_k\in R(X_i,X_j)\}\le K
 \end{aligned}
 \right\}.
\]

We claim that whenever $\mathcal A$ occurs, there is an ordered pair $(i,j)$
whose event $E_{ij}$ serves as a witness. In fact, choose any
$x\in S$ and a point $y\in S$ farthest from $x$.
Then $S\subseteq B(x,s)$.
Thus $R(x,y)$ contains no point of $S$.
Every sample point in $R(x,y)$ is adjacent to $x$ or $y$, so it is
an external neighbor of $S$. There are at most $K$ such points.
Consequently, the labels of $x,y$ satisfy $E_{ij}$, and
\[
 \mathcal A\subseteq\bigcup_{i\ne j}E_{ij}.
\]
This reduces the problem to estimating the probability of each
ordered-pair event.

\begin{figure}[htbp]
\centering
\begin{tikzpicture}[
 scale=1.6,font=\small,line cap=round,line join=round]
 \fill[blue!16] (-.28,0) circle (1.8);
 \fill[blue!16] ( .28,0) circle (1.8);
 \fill[white] (-.28,0) circle (.56);

 \draw[blue!65!black] (-.28,0) circle (1.8);
 \draw[blue!65!black] ( .28,0) circle (1.8);
 \draw[black!65,densely dashed] (-.28,0) circle (.56);

 \fill (-.28,0) circle (.025) node[below left] {$x$};
 \fill ( .28,0) circle (.025) node[below right] {$y$};
 \fill (-.14,.14) circle (.022);
 \fill (-.01,.03) circle (.022);
 \fill (-.09,-.16) circle (.022);
 \node at (-.28,.34) {$S$};

 \node[blue!65!black] at (-1.4,1.15) {$C_x$};
 \node[blue!65!black] at ( 1.4,1.15) {$C_y$};
 \node[blue!65!black] at (0,1.05) {$R(x,y)$};

 \node (core) at (-1.2,-.95) {$B(x,s)$};
 \draw[->,black!65] (core.north east)--(-.67,-.40);
 \draw[<->,black!65] (-.28,-.78)--(.28,-.78)
   node[midway,below] {$s=\dist(x,y)$};
\end{tikzpicture}
\caption{The ordered-pair witness.}
\label{fig:local-tiny-pair}
\end{figure}
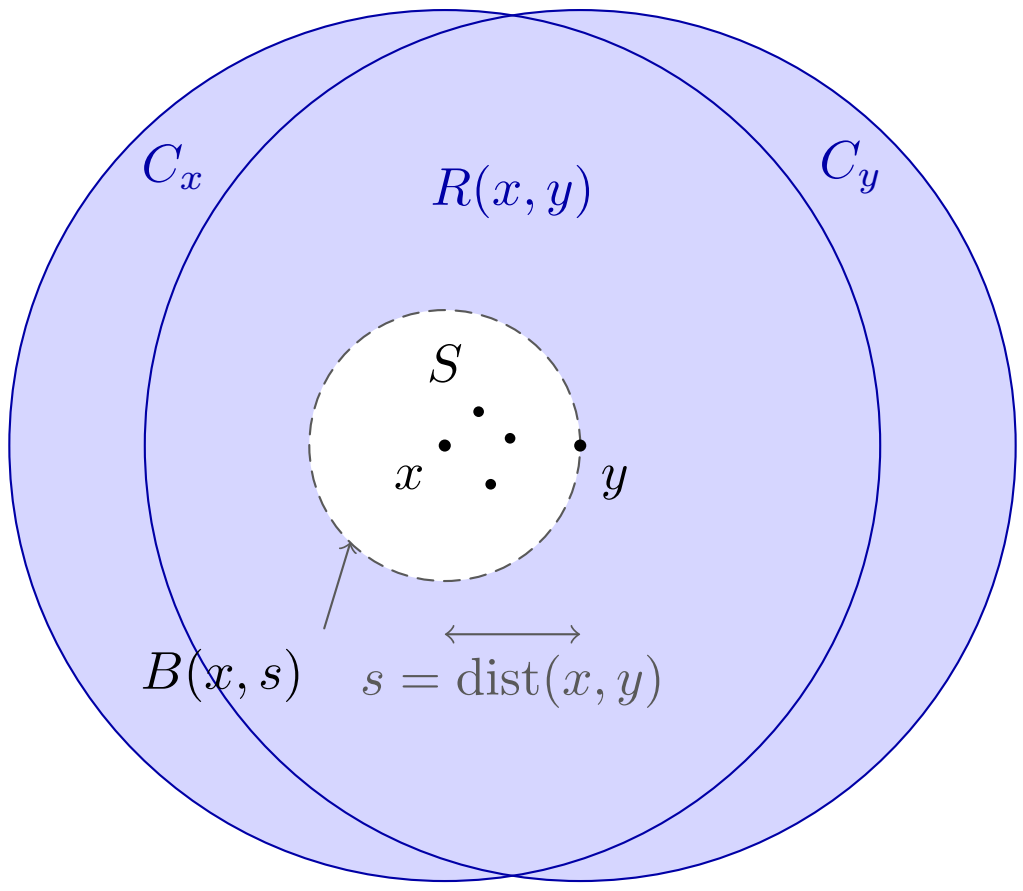

We now fix the two designated locations $x,y$ and put
$u=s/r\in(0,\delta]$. Since $B(x,s)\subseteq C_x$,
\eqref{eq:near-cap} and \eqref{eq:cap-union} give
\[
 2p\ge\mu(R(x,y))
 =\mu(C_x\cup C_y)-F(ur)\ge p(1+cu),
\]
where $c=\frac{1}{\pi}-\frac{1}{5}$.

Conditional on $X_i=x$ and $X_j=y$, the other $n-2$ sample points
are independent and uniform. Their number in $R(x,y)$ therefore
has distribution $\operatorname{Bin}(n-2,\mu(R(x,y)))$.
Applying \eqref{eq:local-binomial} gives
\[
 \Pp(E_{ij}\mid X_i=x,\ X_j=y)
 \le C_K(\log n)^K
       e^{-(n-2)p}e^{-(n-2)cpu}.
\]
For fixed $X_i=x$, the probability that
$\dist(x,X_j)/r\le u$ is $F(ur)$. Integrating over the second
location and applying the union bound over ordered labels yields
\[
 \Pp(\mathcal A)
 \le\sum_{i\ne j}\Pp(E_{ij})
 \le C_K(\log n)^K n^2e^{-(n-2)p}
       \int_0^\delta e^{-(n-2)cpu}\,dF(ur).
\]

Put $a=(n-2)cp\asymp np$. Stieltjes integration by parts and
$F(ur)\le pu^2$ give
\[
\begin{aligned}
 \int_0^\delta e^{-au}\,dF(ur)
 &=e^{-a\delta}F(\delta r)
   +a\int_0^\delta e^{-au}F(ur)\,du\\
 &\le p\delta^2e^{-a\delta}
   +ap\int_0^\delta u^2e^{-au}\,du\\
 &\le Cp/a^2.
\end{aligned}
\]
Substituting this estimate and using \eqref{eq:local-baseline},
we obtain
\begin{equation}\label{eq:local-tiny}
 \Pp(\mathcal A)
 \le C_K(\log n)^K\frac{ne^{-np}}{np}
 \le \frac{C_{K,T}}{\log n}.
\end{equation}

\subsection{A three-point region in the lower-dimensional range}

We now assume $d-1\le\sqrt{\log n}$ and consider a $K$-separated
set $S$ with
\[
 2\le|S|\le b,\qquad D:=\diam(S)\ge\delta r.
\]
Our goal is to select three points $a,x,y\in S$ and construct
a region $R=R(a,x,y)$ satisfying
\[
 R\cap S=\varnothing,\qquad
 R\subseteq C_a\cup C_x\cup C_y,\qquad
 \mu(R)\ge(1+1/200)p.
\]
These properties imply that $R$ contains at most $K$ sample points,
all external neighbors of $S$.

The radius bound in Lemma~\ref{lem:cap-estimates} gives
\[
 r\le C(\log n/n)^{1/(d-1)}
 \le C\exp(-\sqrt{\log n}+o(1)).
\]
Set
\[
 \rho=(b+1)r,\qquad r_0=r\sqrt{1-\rho^2}.
\]
Then $D\le\rho=o(1)$, hence \(S\) contains no antipodal points.

The construction proceeds as follows.

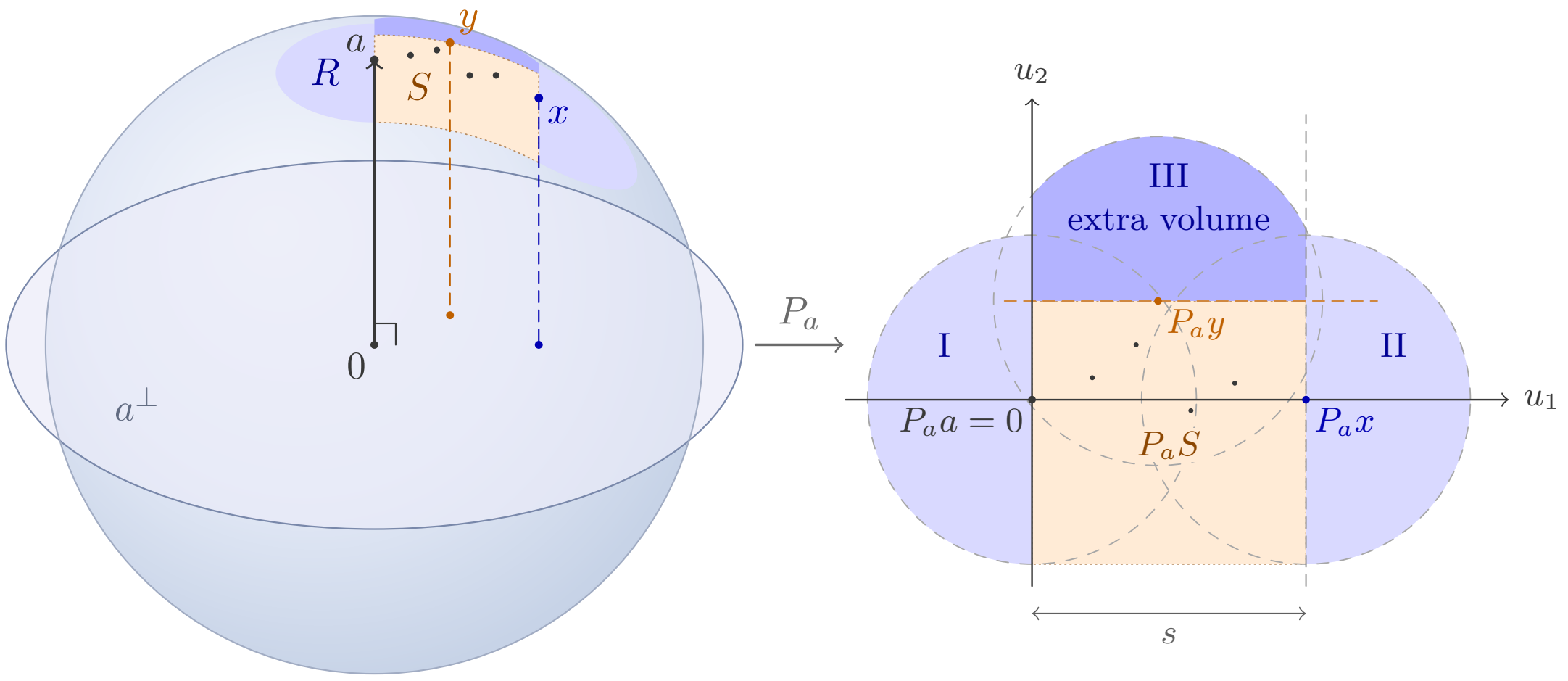
\begin{figure}[tbp]
\centering
\begin{tikzpicture}[font=\small,line cap=round,line join=round]

 \definecolor{capSphereLight}{HTML}{F5F7FC}
 \definecolor{capSphereMid}{HTML}{E3E9F6}
 \definecolor{capSphereShade}{HTML}{BECCE2}
 \definecolor{capPlane}{HTML}{EBECF9}
 \definecolor{capPlaneEdge}{HTML}{7988AA}
 \colorlet{capRegion}{blue!15}
 \colorlet{capExtra}{blue!30}
 \colorlet{capRegionInk}{blue!58!black}
 \colorlet{capSet}{orange!16}
 \colorlet{capSetInk}{orange!55!black}
 \colorlet{capX}{blue!70!black}
 \colorlet{capY}{orange!75!black}
 \colorlet{capInk}{black!78}
 \colorlet{capGuide}{black!43}
 \pgfdeclareradialshading{threePointSphere}{\pgfpoint{-8bp}{9bp}}{
   color(0bp)=(capSphereLight);
   color(17bp)=(capSphereMid);
   color(33bp)=(capSphereShade);
   color(50bp)=(capSphereShade)}

 \pgfmathsetmacro{\endpointX}{0.50}
 \pgfmathsetmacro{\thirdX}{0.23}
 \pgfmathsetmacro{\thirdY}{0.18}
 \pgfmathsetmacro{\ballRadius}{0.30}
 \pgfmathsetmacro{\stripBottom}{-0.30}
 \tikzset{sample patch/.style={fill=capSet,draw=capSetInk!62!white,
   line width=.33pt,densely dotted}}
 \pgfmathsetmacro{\cameraSin}{0.50}
 \pgfmathsetmacro{\cameraCos}{sqrt(3)/2}
 \pgfmathsetmacro{\upperStart}{
   acos((\endpointX-\thirdX)/\ballRadius)}
 \pgfmathsetmacro{\upperEnd}{
   acos(-\thirdX/\ballRadius)}

 \begin{scope}[scale=3,declare function={
   capheight(\cu,\cv)=\cameraSin*\cv+\cameraCos*sqrt(1-\cu*\cu-\cv*\cv);}]
   \shade[shading=threePointSphere] (0,0) circle (1);

   \fill[capPlane,fill opacity=.69] (0,0)
     ellipse[x radius=1.12,y radius=.56];
   \draw[capPlaneEdge,line width=.42pt] (0,0)
     ellipse[x radius=1.12,y radius=.56];
   \node[capPlaneEdge!65!capInk] at (-.72,-.18) {$a^\perp$};

   \coordinate (sphereA) at (0,\cameraCos);
   \coordinate (sphereX) at
     (\endpointX,{\cameraCos*sqrt(1-\endpointX^2)});
   \coordinate (sphereY) at
     (\thirdX,
      {\cameraSin*\thirdY
       +\cameraCos*sqrt(1-\thirdX^2-\thirdY^2)});

   \begin{scope}
     \clip (0,0) circle (1);

     \fill[capRegion]
       plot[domain=90:270,samples=55,variable=\plotAngle]
       ({\ballRadius*cos(\plotAngle)},
        {\cameraSin*\ballRadius*sin(\plotAngle)
         +\cameraCos*sqrt(1-\ballRadius^2)})
       -- cycle;

     \fill[capRegion]
       plot[domain=-90:90,samples=55,variable=\plotAngle]
       ({\endpointX+\ballRadius*cos(\plotAngle)},
        {\cameraSin*\ballRadius*sin(\plotAngle)
         +\cameraCos*sqrt(
           1-(\endpointX+\ballRadius*cos(\plotAngle))^2
            -(\ballRadius*sin(\plotAngle))^2)})
       -- cycle;

     \fill[capExtra]
       plot[domain=\upperStart:\upperEnd,
            samples=65,variable=\plotAngle]
       ({\thirdX+\ballRadius*cos(\plotAngle)},
        {\cameraSin*(\thirdY+\ballRadius*sin(\plotAngle))
         +\cameraCos*sqrt(
           1-(\thirdX+\ballRadius*cos(\plotAngle))^2
            -( \thirdY+\ballRadius*sin(\plotAngle))^2)})
       -- plot[domain=0:\endpointX,
               samples=40,variable=\plotCoordinate]
       (\plotCoordinate,
        {\cameraSin*\thirdY
         +\cameraCos*sqrt(1-\plotCoordinate^2-\thirdY^2)})
       -- cycle;

     \filldraw[sample patch]
       plot[domain=0:\endpointX,samples=45,variable=\cu]
         (\cu,{capheight(\cu,\stripBottom)})
       -- plot[domain=\endpointX:0,samples=45,variable=\cu]
         (\cu,{capheight(\cu,\thirdY)})
       -- cycle;
   \end{scope}

   \draw[capPlaneEdge!70!white,line width=.4pt] (0,0) circle (1);

   \draw[->,capInk,line width=.7pt] (0,0)--(sphereA);
   \draw[capInk] (0,.065)--(.065,.065)--(.065,0);

   \draw[capX,densely dashed]
     (sphereX)--(\endpointX,0);
   \draw[capY,densely dashed]
     (sphereY)--(\thirdX,{\cameraSin*\thirdY});

   \fill[capInk] (0,0) circle (.012)
     node[below left,inner sep=2pt] {$0$};
   \fill[capX] (\endpointX,0) circle (.012);
   \fill[capY]
     (\thirdX,{\cameraSin*\thirdY}) circle (.012);

   \foreach \sampleU/\sampleV in
     {.11/.04,.19/.10,.29/-.02,.37/.03}
   {
     \fill[capInk]
       (\sampleU,
        {\cameraSin*\sampleV
         +\cameraCos*sqrt(1-\sampleU^2-\sampleV^2)})
       circle (.010);
   }

   \fill[capInk] (sphereA) circle (.013)
     node[above left,inner sep=2pt] {$a$};
   \fill[capX] (sphereX) circle (.013)
     node[below right,inner sep=2pt] {$x$};
   \fill[capY] (sphereY) circle (.013)
     node[above right,inner sep=2pt] {$y$};

   \node[capSetInk] at (.135,{capheight(.135,-.13)}) {$S$};
   \node[capRegionInk] at (-.15,.83) {$R$};

   \node[below] at (0,-1.16)
     {(a) The sphere, \(S\), and \(R\)};
 \end{scope}

 \draw[->,capInk!75!white,line width=.65pt]
   (3.48,0)--(4.28,0)
   node[midway,above] {$P_a$};

 \begin{scope}[shift={(6,-.5)},scale=5,font=\footnotesize]

   \filldraw[sample patch]
     (0,\stripBottom) rectangle (\endpointX,\thirdY);

   \begin{scope}
     \clip (-.34,-.34) rectangle (0,.52);
     \fill[capRegion] (0,0) circle (\ballRadius);
   \end{scope}

   \begin{scope}
     \clip (\endpointX,-.34) rectangle (.84,.52);
     \fill[capRegion]
       (\endpointX,0) circle (\ballRadius);
   \end{scope}

   \begin{scope}
     \clip (0,\thirdY) rectangle (\endpointX,.52);
     \fill[capExtra]
       (\thirdX,\thirdY) circle (\ballRadius);
   \end{scope}

   \draw[capGuide!80!white,line width=.35pt,dashed] (0,0) circle (\ballRadius);
   \draw[capGuide!80!white,line width=.35pt,dashed]
     (\endpointX,0) circle (\ballRadius);
   \draw[capGuide!80!white,line width=.35pt,dashed]
     (\thirdX,\thirdY) circle (\ballRadius);

   \draw[capGuide,line width=.38pt,densely dashed]
     (\endpointX,-.34)--(\endpointX,.52);
   \draw[capY!75!white,densely dashed]
     (-.05,\thirdY)--(.63,\thirdY);

   \draw[->,capInk,line width=.45pt] (-.34,0)--(.87,0)
     node[right] {$u_1$};
   \draw[->,capInk,line width=.45pt] (0,-.34)--(0,.55)
     node[above] {$u_2$};

   \foreach \sampleU/\sampleV in
     {.11/.04,.19/.10,.29/-.02,.37/.03}
   {
     \fill[capInk] (\sampleU,\sampleV) circle (.0045);
   }

   \fill[capInk] (0,0) circle (.007)
     node[below left,inner sep=2pt] {$P_a a=0$};
   \fill[capX] (\endpointX,0) circle (.007)
     node[below right,inner sep=2pt] {$P_a x$};
   \fill[capY] (\thirdX,\thirdY) circle (.007)
     node[below right,inner sep=2pt] {$P_a y$};

   \node[capSetInk,fill=capSet,inner sep=.8pt] at (.25,-.083) {$P_aS$};
   \node[capRegionInk] at (-.16,.10) {$\mathrm{I}$};
   \node[capRegionInk] at ({\endpointX+.16},.10)
     {$\mathrm{II}$};
   \node[capRegionInk,align=center] at (.25,.37)
     {$\mathrm{III}$\\extra volume};

   \draw[<->,capInk!80!white] (0,-.39)--(\endpointX,-.39)
     node[midway,below] {$s$};

   \node[below] at (.25,-.59)
     {(b) Coordinates, \(P_aS\), and \(\widetilde R=P_aR\)};
 \end{scope}

\end{tikzpicture}
\caption{The three-point construction.}
\label{fig:local-projection}
\label{fig:slab-certificate}
\end{figure}

\begin{enumerate}
\item
Choose diameter endpoints $a,x\in S$ and project onto $a^\perp$ by
\[
 P_a z=z-(a\cdot z)a.
\]
The point $a$ projects to the origin. Choose the first coordinate
axis along $P_a x$, so that
\[
 P_a x=(s,0,\ldots,0),\qquad s=\|P_a x\|.
\]
Choose the second coordinate direction by a fixed measurable rule
depending only on $a,x$.

\item
Choose $y\in S$ whose projection $P_a y$ has largest second
coordinate.

\item
In $a^\perp$, write
\[
\begin{aligned}
 \mathrm{I}
 &=B^{d-1}(0,r_0)\cap\{u:u_1<0\},\\
 \mathrm{II}
 &=B^{d-1}(P_a x,r_0)\cap\{u:u_1>s\},\\
 \mathrm{III}
 &=B^{d-1}(P_a y,r_0)\cap
   \{u:0<u_1<s,\ u_2>(P_a y)_2\}.
\end{aligned}
\]
as marked in
Figure~\ref{fig:local-projection}(b). Define
\[
\widetilde R:=\mathrm{I}\cup\mathrm{II}\cup\mathrm{III}.
\]
The first two pieces are outward half-balls at the endpoints.
The third lies above $P_a y$ and between the two vertical boundaries.

\item
Lift $\widetilde R$ from the unit ball in $a^\perp$ to the hemisphere
centered at $a$ via the map
\[
\begin{aligned}
 \Phi_a \colon B^{d-1}(0,1)
 &\longrightarrow \{z\in\mathbb S^{d-1}:a\cdot z\ge0\},\\
 u &\longmapsto u+\sqrt{1-\|u\|^2}\,a.
\end{aligned}
\]
Set $R=\Phi_a(\widetilde R)$.
\end{enumerate}

Figure~\ref{fig:local-projection} shows the process of this construction. The coordinate choice and the region depend only on $a,x,y$ and
the fixed parameters.  We next verify its
geometric properties and then estimate its area.

\begin{lemma}[Location of the three-point region]\label{lem:local-region}
For the construction above, $\widetilde R$ avoids $P_aS$.
Its image $R$ is well defined and satisfies
\[
 R\cap S=\varnothing,\qquad
 R\subseteq C_a\cup C_x\cup C_y,\qquad
 \mu(R)\le3p.
\]
\end{lemma}

\begin{proof}
Connectivity and $|S|\le b$ give $D\le br$.
Since $a,x$ are diameter endpoints, for every $z\in S$,
\[
 a\cdot z,\ x\cdot z\ge a\cdot x=1-D^2/2>0.
\]
Denote
\[
 s:=\|P_a x\|=D\sqrt{1-D^2/4}.
\]
In the coordinates chosen above,
\begin{equation}\label{eq:local-spherical-slab}
 0\le\frac{x\cdot z-(a\cdot x)(a\cdot z)}s=(P_a z)_1
 \le\|P_a z\|\le s.
\end{equation}

Thus $P_aS$ lies between the two vertical boundaries.
By the choice of $y$, it also lies below the horizontal boundary
through $P_a y$. Each of the three pieces defining $\widetilde R$
therefore avoids $P_aS$.

The three
radius-$r_0$ balls lie in $\{\|u\|\le br+r_0\}\subseteq\{\|u\|\le\rho\}$,
where $\rho<1$. Thus the map $\Phi_a$ is well-defined on
$\widetilde R$.
Differentiating $\Phi_a$ gives the distance bound
\[
 \|\Phi_a(u)-\Phi_a(u')\|
 \le\frac{\|u-u'\|}{\sqrt{1-\rho^2}}
 \qquad(\|u\|,\|u'\|\le\rho).
\]
Since $r_0/\sqrt{1-\rho^2}=r$, each planar ball maps into the
radius-$r$ cap about its corresponding point. Consequently,
\[
 R=\Phi_a(\widetilde R)\subseteq C_a\cup C_x\cup C_y.
\]
The inverse relation  gives
$R\cap S=\varnothing$ and it's easy to check
$\mu(R)\le3p$.
\end{proof}

\begin{lemma}[Fixed area gain]\label{lem:slab-certificate}
For all sufficiently large $n$, uniformly in the stated dimension
range, the construction above satisfies
\[
 |\widetilde R|\ge1.01\,|B^{d-1}(0,r_0)|,
 \qquad
 \mu(R)\ge(1+1/200)p.
\]
\end{lemma}

\begin{proof}
The cuts defining Regions $\mathrm{I}$ and $\mathrm{II}$ pass
through their respective ball centers, so
\[
 |\mathrm{I}|+|\mathrm{II}|
 =\frac12|B^{d-1}(0,r_0)|+\frac12|B^{d-1}(0,r_0)|
 =|B^{d-1}(0,r_0)|.
\]
It remains to estimate the additional volume supplied by
Region $\mathrm{III}$.

Fix $u_3,\ldots,u_{d-1}$ so that the corresponding section of
$B^{d-1}(P_a y,r_0)$ has positive area. This section is a disk of radius
$0<h\le r_0$ in the $(u_1,u_2)$ coordinates. We show that
Region~$\mathrm{III}$ occupies more than $1\%$ of its area.

Note that the center $P_a y$ has first coordinate in $[0,s]$.
On at least one side of that coordinate, the interval $(0,s)$
extends a distance of at least $s/2$. Moreover,
\[
 s=D\sqrt{1-D^2/4}
 \ge\delta r\sqrt{1-\rho^2}
 =r_0/10.
\]

Thus a one-sided interval of length $r_0/20$, measured from the
ball center, is available inside the strip $0<u_1<s$.

The third piece retains the part of this disk above its horizontal
diameter and inside the strip $0<u_1<s$. Hence
\[
\begin{aligned}
 \frac{\text{retained sectional area}}{\pi h^2}
 &\ge
 \frac1{\pi h^2}
 \int_0^{\min\{r_0/20,h\}}\sqrt{h^2-t^2}\,dt\\
 &=\frac1\pi
 \int_0^{\min\{r_0/(20h),1\}}\sqrt{1-t^2}\,dt\\
 &\ge\frac1\pi\int_0^{1/20}\sqrt{1-t^2}\,dt
 >0.01.
\end{aligned}
\]
Thus every two-dimensional section retains more than $1\%$ of
its disk area. Integrating over the remaining coordinates shows
that the third piece has volume at least
\(0.01\,|B^{d-1}(0,r_0)|.
\)
Adding the two half-balls  $\mathrm{I}$ and $\mathrm{II}$ yields
\[
 |\widetilde R|\ge1.01\,|B^{d-1}(0,r_0)|.
\]

The area factor of the lift $\Phi_a$ is
$(1-\|u\|^2)^{-1/2}\ge1$. Dividing by the sphere's area and using \eqref{eq:geometry_for_w_d}
gives
\[
 \mu(R)
 \ge \frac{|\widetilde R|}{|\mathbb S^{d-1}|}
 \ge 1.01\,\frac{|B^{d-1}(0,r_0)|}{|\mathbb S^{d-1}|}
 = 1.01\,w_d r_0^{d-1}.
\]
The density \eqref{eq:chord-density} gives $p\le w_dr^{d-1}$.
Therefore
\[
 \mu(R)
 \ge1.01(1-\rho^2)^{(d-1)/2}p
 \ge(1+1/200)p
\]
for sufficiently large $n$, since $(d-1)\rho^2=o(1)$ uniformly.
\end{proof}

Set \(g=1/200\). Let $\mathcal B$ be the event that there exists a $K$-separated set
$S$ with $2\le|S|\le b$ and $\diam(S)\ge\delta r$.
For $i\ne j$ and $k\in[n]$, let $R_{ijk}$ be the region constructed
from $X_i,X_j,X_k$, with $R_{ijk}=\varnothing$ if the construction
is undefined. Define
\[
 E_{ijk}=
 \left\{
 \begin{aligned}
 &X_j,X_k\in B(X_i,br),\qquad
   (1+g)p\le\mu(R_{ijk})\le3p,\\
 &\#\{\ell\in[n]:X_\ell\in R_{ijk}\}\le K
 \end{aligned}
 \right\}.
\]
Every set $S$ witnessing $\mathcal B$ supplies such a triple.
Consequently,
\[
 \mathcal B\subseteq\bigcup_{i\ne j}\bigcup_{k=1}^n E_{ijk}.
\]

For fixed $i\ne j$ and $k$, condition on $X_i,X_j,X_k$.
The region $R_{ijk}$ is then fixed, and the remaining sample points
are independent and uniform. Applying \eqref{eq:local-binomial} to
these points and averaging over $X_i,X_j,X_k$ gives
\[
\begin{aligned}
 \Pp(\mathcal B)
 &\le\sum_{i\ne j}\sum_{k=1}^n\Pp(E_{ijk})\\
 &\le C_K(\log n)^K e^{-(n-3)(1+g)p}
       \left[n^3F(br)^2+2n^2F(br)\right]\\
 &\le C_K(\log n)^K n[nF(br)]^2e^{-(n-3)(1+g)p}.
\end{aligned}
\]
The two terms account for three distinct labels and for $k=i$ or
$k=j$, respectively. The last step uses $nF(br)\ge np\to\infty$.

By \eqref{eq:near-cap}, with $\lambda=b$,
\begin{equation}\label{eq:local-flat}
 nF(br)\le np(\sqrt2\,b)^{d-1}
 =\exp\{O_{K,T}(\sqrt{\log n}\log\log n)\}.
\end{equation}
Together with \eqref{eq:local-baseline}, this yields
\[
 \Pp(\mathcal B)
 \le\exp\{-g\log n+O_{K,T}(\sqrt{\log n}\log\log n)\}
 \le n^{-g/2}.
\]
Combining this bound with \eqref{eq:local-tiny} proves
Proposition~\ref{prop:local} throughout the lower-dimensional range.

\subsection{Dimensions \texorpdfstring{$d\ge\sqrt{\log n}$}{d at least the square root of log n}}

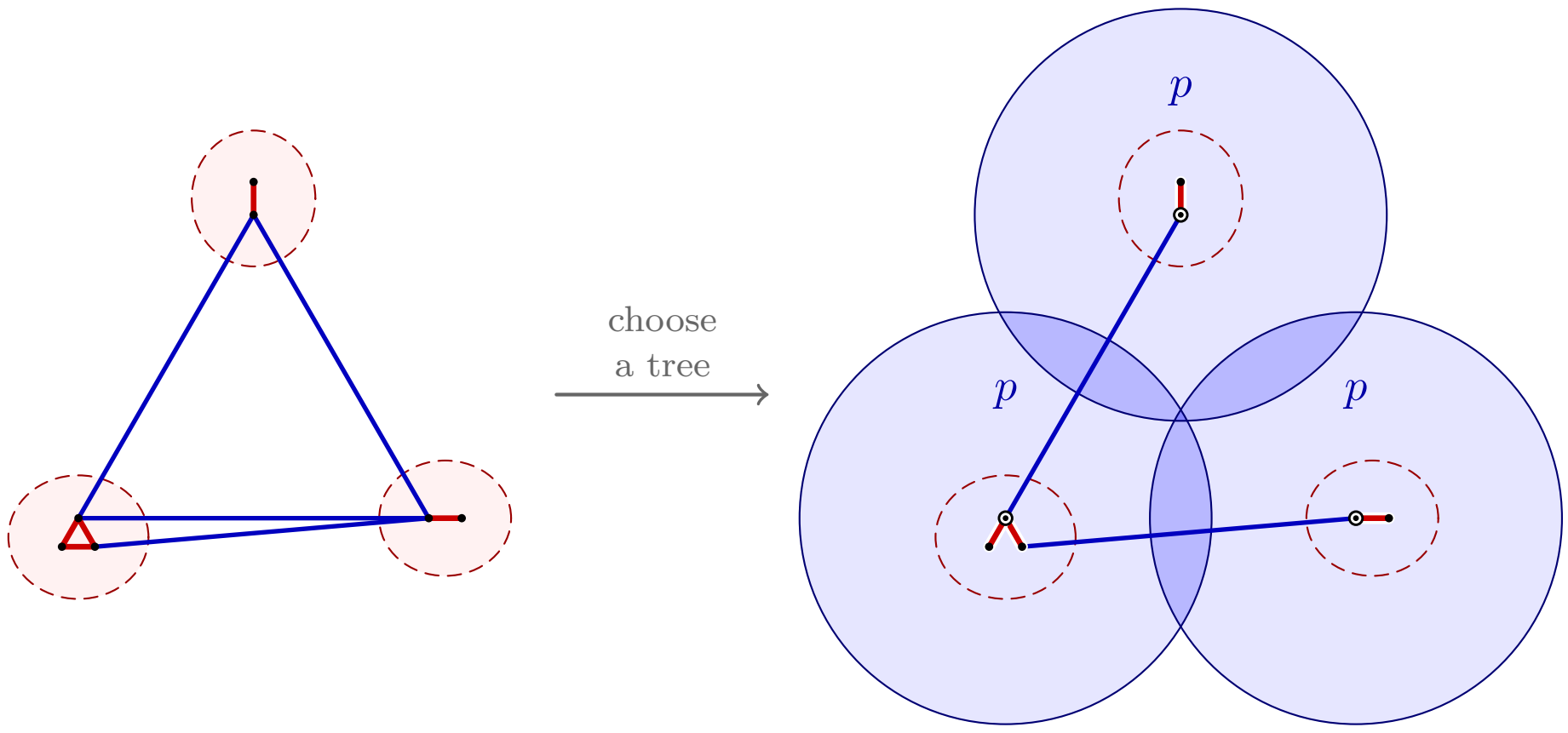
\begin{figure}[tb]
\centering
\begin{tikzpicture}[
 scale=1.60,font=\small,line cap=round,line join=round,
 red edge/.style={red!80!black,line width=1.25pt},
 tree red edge/.style={
   red edge,preaction={draw=white,line width=2.6pt}},
 blue edge/.style={blue!75!black,line width=1pt}]

 \pgfmathsetmacro{\clusterSide}{1.70}
 \pgfmathsetmacro{\clusterHeight}{sqrt(3)*\clusterSide/2}
 \pgfmathsetmacro{\topCenterX}{\clusterSide/2}
 \pgfmathsetmacro{\clusterOffset}{.16}
 \pgfmathsetmacro{\lowerLeafY}{-sqrt(3)*\clusterOffset/2}
 \pgfmathsetmacro{\capRadius}{1}

 \begin{scope}
   \coordinate (L1) at (0,0);
   \coordinate (L2) at
     ({\clusterOffset/2},\lowerLeafY);
   \coordinate (L3) at
     ({-\clusterOffset/2},\lowerLeafY);
   \coordinate (L4) at (\clusterSide,0);
   \coordinate (L5) at
     ({\clusterSide+\clusterOffset},0);
   \coordinate (L6) at (\topCenterX,\clusterHeight);
   \coordinate (L7) at
     (\topCenterX,{\clusterHeight+\clusterOffset});

   \fill[red!5] (0,-.092)
     ellipse[x radius=.34,y radius=.30];
   \fill[red!5] ({\clusterSide+.08},0)
     ellipse[x radius=.32,y radius=.28];
   \fill[red!5] (\topCenterX,{\clusterHeight+.08})
     ellipse[x radius=.30,y radius=.33];

   \draw[red!60!black,densely dashed] (0,-.092)
     ellipse[x radius=.34,y radius=.30];
   \draw[red!60!black,densely dashed]
     ({\clusterSide+.08},0)
     ellipse[x radius=.32,y radius=.28];
   \draw[red!60!black,densely dashed]
     (\topCenterX,{\clusterHeight+.08})
     ellipse[x radius=.30,y radius=.33];

   \draw[blue edge] (L1)--(L4);
   \draw[blue edge] (L2)--(L4);
   \draw[blue edge] (L1)--(L6);
   \draw[blue edge] (L4)--(L6);

   \draw[red edge] (L1)--(L2)--(L3)--cycle;
   \draw[red edge] (L4)--(L5);
   \draw[red edge] (L6)--(L7);

   \foreach \vertex in {1,...,7}
     \fill (L\vertex) circle (.020);

   \node[below] at (\topCenterX,-1.65)
     {(a) Connected configuration};
 \end{scope}

 \draw[->,thick,black!60] (2.32,.60)--(3.35,.60)
   node[midway,above,align=center,font=\scriptsize]
   {choose\\a tree};

 \begin{scope}[shift={(4.5,0)}]
   \coordinate (T1) at (0,0);
   \coordinate (T2) at
     ({\clusterOffset/2},\lowerLeafY);
   \coordinate (T3) at
     ({-\clusterOffset/2},\lowerLeafY);
   \coordinate (T4) at (\clusterSide,0);
   \coordinate (T5) at
     ({\clusterSide+\clusterOffset},0);
   \coordinate (T6) at (\topCenterX,\clusterHeight);
   \coordinate (T7) at
     (\topCenterX,{\clusterHeight+\clusterOffset});

   \foreach \vertex in {1,4,6}
     \fill[blue!10] (T\vertex) circle (\capRadius);

   \begin{scope}
     \clip (T1) circle (\capRadius);
     \fill[blue!28] (T4) circle (\capRadius);
     \fill[blue!28] (T6) circle (\capRadius);
   \end{scope}
   \begin{scope}
     \clip (T4) circle (\capRadius);
     \fill[blue!28] (T6) circle (\capRadius);
   \end{scope}

   \begin{scope}
     \clip (T1) circle (\capRadius);
     \clip (T4) circle (\capRadius);
     \fill[blue!48] (T6) circle (\capRadius);
   \end{scope}

   \foreach \vertex in {1,4,6}
     \draw[blue!45!black]
       (T\vertex) circle (\capRadius);

   \node[blue!65!black] at (0,.60) {$p$};
   \node[blue!65!black] at (\clusterSide,.60) {$p$};
   \node[blue!65!black] at
     (\topCenterX,{\clusterHeight+.60}) {$p$};

   \draw[red!60!black,densely dashed] (0,-.092)
     ellipse[x radius=.34,y radius=.30];
   \draw[red!60!black,densely dashed]
     ({\clusterSide+.08},0)
     ellipse[x radius=.32,y radius=.28];
   \draw[red!60!black,densely dashed]
     (\topCenterX,{\clusterHeight+.08})
     ellipse[x radius=.30,y radius=.33];

   \draw[blue edge] (T2)--(T4);
   \draw[blue edge] (T1)--(T6);

   \draw[tree red edge] (T1)--(T2);
   \draw[tree red edge] (T1)--(T3);
   \draw[tree red edge] (T4)--(T5);
   \draw[tree red edge] (T6)--(T7);

   \foreach \vertex in {1,...,7}
     \fill (T\vertex) circle (.020);

   \foreach \vertex in {1,4,6}
   {
     \draw[black,fill=white,line width=.55pt]
       (T\vertex) circle (.032);
     \fill (T\vertex) circle (.014);
   }

   \node[below] at (\topCenterX,-1.65)
     {(b) Red--blue spanning tree};
 \end{scope}

\end{tikzpicture}
\caption{Reduction to a red--blue spanning tree.}
\label{fig:local-typed-tree}
\end{figure}

Fix a candidate vertex set $S=\{1,\ldots,s\}$ such that $2\le s\le b$. We will encode this set by a spanning tree with two edge types.
Color an edge red when its length is at most $\delta r$, and blue
when its length lies in $(\delta r,r]$.  Retain a spanning tree inside each red component
and then choose blue edges connecting these components into a tree. Figure~\ref{fig:local-typed-tree} illustrates this reduction.

For $1\le m\le s$, define
\[
 \mathcal C_m:=
 \left\{
 \begin{aligned}
 &G_d(n,p)[S]\text{ is connected},\\
 &\text{its red subgraph has exactly }m\text{ components}
 \end{aligned}
 \right\}.
\]

We will estimate the external-neighbor probability conditioning on \(\mathcal C_m\),  and then sum over $m$.

\begin{enumerate}

\item \emph{External-neighbor probability conditional on $\mathcal C_m$.}

Recall that \eqref{eq:near-cap} and \eqref{eq:far-overlap} give absolute constants $a,C_1>0$
such that
\[
 p_\delta\le C_1e^{-d}p,
 \qquad
 \dist(x,y)>\delta r
 \ \Longrightarrow\
 \mu(C_x\cap C_y)\le C_1e^{-a(d\wedge\log n)}p.
\]
Write
\[
 U=\mu\left(\bigcup_{i=1}^s C_{X_i}\right).
\]
On $\mathcal C_m$, choose one vertex from each of the $m$ red
components. The chosen vertices are pairwise more than $\delta r$ apart,
so \eqref{eq:far-overlap} applies to each pair of their caps.

Applying
Bonferroni's inequality to these $m$ chosen caps gives
\[
 U\ge mp-\binom{m}{2}C_1e^{-a(d\wedge\log n)}p.
\]
Consequently, on $\mathcal C_m$,
\[
\begin{aligned}
 -(n-s)U
 &\le-nmp+smp+C_1nm^2e^{-a(d\wedge\log n)}p\\
 &\le-nmp+C(\log n)^{17}
       \left(\frac1n+e^{-a\sqrt{\log n}}\right)\\
 &=-nmp+o(1).
\end{aligned}
\]
The error is uniform over $m\le s\le b$ and
$d\ge\sqrt{\log n}$.

Conditional on $X_1,\ldots,X_s$, the number of external neighbors
has distribution $\operatorname{Bin}(n-s,U)$.
Applying \eqref{eq:local-binomial}, using $U\le sp$, and then
averaging over the locations within $\mathcal C_m$, we obtain
\begin{equation}\label{eq:local-high-tail}
 \Pp\{|\partial S|\le K\mid\mathcal C_m\}\le C_K(s\log n)^K e^{-nmp}.
\end{equation}

\item \emph{Probability of $\mathcal C_m$.}

On $\mathcal C_m$, we can select a spanning tree with $s-m$ red
edges and $m-1$ blue edges.
Cayley's formula gives $s^{s-2}$ labeled trees, and there are
$\binom{s-1}{m-1}$ choices of the blue edges.
For each choice, successive leaf integration contributes
$p_\delta$ per red edge and $p-p_\delta$ per blue edge.
The union bound therefore gives
\[
 \Pp(\mathcal C_m)
\le s^{s-2}\binom{s-1}{m-1}
       p_\delta^{s-m}p^{m-1}.
\]

\item \emph{Summation over $m$ and $s$.}

Let $B_s$ count $K$-separated sets of order $s$.
Note that
\[
 n(p_\delta+pe^{-np})
 =O_{K,T}\bigl((\log n)^{-1}\bigr)
\]
We sum over $m$ by the binomial theorem.
\[
\begin{aligned}
 \sum_{s=2}^b\E B_s
 &=\sum_{s=2}^b\binom ns
   \sum_{m=1}^s
   \Pp\{|\partial S|\le K\mid\mathcal C_m\}
   \Pp(\mathcal C_m)\\
 &\le C_K(\log n)^Ke^{-np}
   \sum_{s=2}^b\binom ns s^{s-2+K}
   \sum_{m=1}^s\binom{s-1}{m-1}
   p_\delta^{s-m}(pe^{-np})^{m-1}\\
 &=C_K(\log n)^Ke^{-np}
   \sum_{s=2}^b\binom ns s^{s-2+K}
   (p_\delta+pe^{-np})^{s-1}\\
 &\le
   C_K(\log n)^Kn^2e^{-np}(p_\delta+pe^{-np})
   \sum_{s=2}^b s^{K-2}
   \bigl[en(p_\delta+pe^{-np})\bigr]^{s-2}\\
 &\le C_{K,T}\,n(p_\delta+pe^{-np})\\
 &=O_{K,T}\bigl((\log n)^{-1}\bigr)
\end{aligned}
\]

\end{enumerate}

Markov's inequality now gives
\[
 \Pp\{\text{there exists a $K$-separated set of order }2,\ldots,b\}
 \le\sum_{s=2}^b\E B_s
 =O_{K,T}\bigl((\log n)^{-1}\bigr)
\]
throughout this dimension range.
The ranges $d-1\le\sqrt{\log n}$ and $d\ge\sqrt{\log n}$
together cover all dimensions, completing the proof of
Proposition~\ref{prop:local}.

\section{Robust structure from independent bridge points}
\label{sec:structure}

Now we show that, with high probability, after any deletion of at most $K$ vertices,
all remaining nonisolated vertices belong to the same component.

For an integer $q\ge1$, let $\mathcal R_q(p)$ be the event that,
simultaneously for every $Z\subseteq[n]$ with $|Z|\le q-1$,
the graph $G-Z$ has exactly one component of order at least two.

\begin{theorem}[Robust structure]\label{thm:robust}
Fix an integer $K\ge0$ and $0\le T<\infty$.
Uniformly over $d\ge3$ and
\[
 |np-\log n-K\log\log n|\le T,
\]
we have
\[
 \Pp\{\mathcal R_{K+1}(p)^c\}
 =O_{K,T}\bigl((\log n)^{-1}\bigr).
\]
\end{theorem}

Thus, with the stated probability, every deletion of at most $K$
vertices leaves one nontrivial component and some isolated vertices.
If the original graph has minimum degree at least $K+1$, no vertex
can become isolated after such a deletion, and the remaining graph
is therefore connected.

We prove the theorem in four steps.
\begin{enumerate}
\item
Split the sample into a large base and a smaller independent reserve.
Apply Proposition~\ref{prop:local} to show that, after any allowed
base deletion, every base component is either a singleton or large.

\item
Show that every partition of the large base components into two
groups has a sufficiently large region of points adjacent to both
groups.

\item
Use the independent reserve points to merge all large base
components, simultaneously for every allowed deletion.

\item
Show that, after each allowed deletion, fewer than $b$ surviving
vertices lie outside the component containing all the large base
components from step~(3). Any other nontrivial component would then
be a small $K$-separated set in the original graph. Apply
Proposition~\ref{prop:local} to rule this out.
\end{enumerate}

\subsection{The base sample and its large components}
\label{sec:structure-base}

Take the first
\[
 M=\left\lfloor\left(1-\frac1{\log n}\right)n\right\rfloor
\]
sample points as the \emph{base set} $A$, and call the remaining points
the \emph{reserve}.
The base graph \(G[A]\) has distribution $G_d(M,p)$, and it's easy to check
\begin{equation}\label{eq:base-window}
 Mp-\log M-K\log\log M=t-1+o(1)=O_{K,T}(1).
\end{equation}

Set
\(
 b=\lceil(\log M)^8\rceil.
\)
Applying Proposition~\ref{prop:local} with sample size $M$ gives
\begin{equation}\label{eq:base-local-exclusion}
 \Pp\left\{
 \begin{array}{c}
 G[A]\text{ contains a $K$-separated set}\\
 \text{of order }2,\ldots,b
 \end{array}
 \right\}
 =O_{K,T}\bigl((\log n)^{-1}\bigr).
\end{equation}

Let
\[
 Y=\{x\in A:\deg_{G[A]}(x)\le K\}.
\]
The degree of a fixed base vertex has distribution
$\operatorname{Bin}(M-1,p)$. Hence
\begin{equation}\label{eq:base-low-degree}
 \E|Y|
 \le C_KM(\log n)^Ke^{-Mp}
 =O_{K,T}(1),
 \qquad
 \Pp\{|Y|>\log n\}
 =O_{K,T}\bigl((\log n)^{-1}\bigr).
\end{equation}

Consequently, with probability $1-O_{K,T}((\log n)^{-1})$,
we have $|Y|\le\log n$, and 
$G[A]$ contains no $K$-separated set $S$ with $2\le |S|\le b$.

\subsection{A bridge between any two groups of large components}
\label{sec:structure-bridges}

For a set $V$ of base vertices,  write
\[
 U_V=\bigcup_{x\in V} C_x.
\]
Given the base, a single reserve point merges two distinct components
$S$ and $T$ with probability $\mu(U_S\cap U_T)$.

\begin{lemma}[Uniform bridge probability]\label{lem:bridge}
There is an event $\mathcal B$ such that
\[
 \Pp(\mathcal B^c)
 =O_{K,T}\bigl((\log n)^{-1}\bigr).
\]
On $\mathcal B$, we have $|Y|\le\log n$, and every deletion
of at most $K$ base vertices leaves only components of order
greater than $b$ and isolated vertices belonging to $Y$.
After any such deletion, whenever all the components of order
greater than $b$ are divided into two nonempty groups with
vertex unions $S,T$, we have
\begin{equation}\label{eq:uniform-bridge}
 \mu(U_S\cap U_T)
 \ge \frac{c(\log n)^{5/2}}{n},
\end{equation}
where $c>0$ is an absolute constant and $n$ is sufficiently large.
\end{lemma}

\begin{proof}
Fix a deletion of at most $K$ base vertices and a partition of the
remaining large components into two nonempty groups with vertex
unions $S,T$. After exchanging the groups if necessary, let
\[
 B=U_S\setminus U_T
\]
be the smaller of the two regions covered by just one group.
Thus $\mu(B)\le 1/2$. We estimate the bridge region $U_S\cap U_T$
by enlarging $B$.

For $0<h<r$, let $B_h$ be the closed chordal $h$-neighborhood of $B$
as in \eqref{def:B_h}. We claim that
\begin{equation}\label{eq:bridge-layer}
 B_h\setminus B
 \subseteq (U_S\cap U_T)
 \,\cup\, \{z:\dist(z,S\cup T)\ge r-h\}.
\end{equation}
In fact, if $z\in B_h\setminus B$ belongs to $U_S$, then it also belongs
to $U_T$.
Otherwise, $\dist(z,S)>r$, and  by the triangle inequality,
$\dist(z,T)\ge r-h$.

Writing $H_V=\{z:\dist(z,V)\ge r-h\}$, \eqref{eq:bridge-layer} yields
\begin{equation}\label{eq:bridge-growth}
 \mu(U_S\cap U_T)
 \ge \mu(B_h)-\mu(B)-\mu(H_{S\cup T}).
\end{equation}

We first choose $h$ so that the last term is small. Take
\[
 h=\frac{r}{\sqrt{d-1}\,(\log n)^{3/2}}.
\]
By \eqref{eq:d-tau^2}, we have $d\tau^2\asymp d\wedge\log n$.
Applying the lower bound in \eqref{eq:mills} at $u=\tau$ and using
$p=a_dJ(\tau),1-\tau^2\asymp r^2$, we obtain
\[\begin{aligned}
 \frac{rF'(r)}{p}
 &= \frac{r^2(1-\tau^2)^{(d-3)/2}}{J(\tau)}\\
 &\le \frac{(d-1)r^2}{1-\tau^2}
       \sqrt{\tau^2+\frac{2(1-\tau^2)}{d+1}}\\
 &\lesssim \frac{d\tau r^2}{1-\tau^2}\lesssim \sqrt{d\log n},
\end{aligned}
\]
By \eqref{eq:chord-density},
$F'(s)=a_d s\bigl(s\sqrt{1-s^2/4}\bigr)^{d-3}$ is nondecreasing
on $[0,r]$, since $d\ge3$ and $r\le\sqrt2$. Hence the outer
layer of a single cap has mass
\[
 p-F(r-h)=\int_{r-h}^{r}F'(s)\,ds
 \le hF'(r)
 \lesssim \frac{p}{\log n}
 =O(1/n).
\]
For a fixed $z$, membership in $H_A$ means that all $M$ base
points are at distance at least $r-h$ from $z$. By independence,
\[
 \E\mu(H_A)
 =(1-F(r-h))^M
 \le \exp\{-Mp+O(1)\}
 \le \frac{C_{K,T}}{n(\log n)^K}.
\]
Markov's inequality now gives
\begin{equation}\label{eq:base-uncovered}
 \Pp\left\{\mu(H_A)>\frac{\log n}{n}\right\}
 =O_{K,T}\bigl((\log n)^{-1}\bigr).
\end{equation}
If $G[A]$ contains no $K$-separated set of order $2,\ldots,b$,
then $S\cup T$ omits at most $K+|Y|$ base vertices. Hence with probability
$1-O_{K,T}((\log n)^{-1})$ we have, for every allowed deletion
and partition,
\begin{equation}\label{eq:base-deletion-uncovered}
 \mu(H_{S\cup T})
 \le \mu(H_A)+(K+|Y|)p
 \le \frac{\log n}{n}+(K+\log n)p
 \le \frac{C_{K,T}(\log n)^2}{n}.
\end{equation}

Now we bound the enlargement gain in \eqref{eq:bridge-growth}. Fix a base
label $i$ and condition on $X_i=x$. For $j\ne i$, the variables
$W_j=q(x,X_j)/p$ are independent and
have mean $p$ by Lemma~\ref{lem:overlap}. Since
$e^u\le 1+(e-1)u$ on $[0,1]$, exponential Markov's inequality
gives
\[
 \Pp\left\{\sum_{j\ne i}W_j>C\log n\,\middle|\,X_i=x\right\}
 \le \exp\{-C\log n+(e-1)(M-1)p\}
 \le n^{-32}
\]
for a sufficiently large absolute $C$. Taking a union bound over
$i$ and including the diagonal term $q(x,x)=p$, we obtain
\[
 \sum_{y\in A}q(x,y)\le Cp\log n
 \qquad\text{for every }x\in A
\]
with probability at least $1-n^{-30}$. Whenever this bound holds,
Cauchy--Schwarz gives, simultaneously for every nonempty $V\subseteq A$,
\begin{equation}\label{eq:base-cap-mass}
 \mu(U_V)
 \ge
 \frac{\bigl(\int\sum_{x\in V}\mathbf{1}_{C_x}\,d\mu\bigr)^2}
      {\int\bigl(\sum_{x\in V}\mathbf{1}_{C_x}\bigr)^2\,d\mu}
 =
 \frac{(|V|p)^2}{\sum_{x,y\in V}q(x,y)}
 \ge \frac{|V|}{C_0n},
\end{equation}
where the last inequality uses $p\asymp(\log n)/n$.

Let $\mathcal B$ be the event that \eqref{eq:base-cap-mass} holds
for every nonempty $V\subseteq A$, that $\mu(H_A)\le(\log n)/n$,
that $|Y|\le\log n$, and that the base has no $K$-separated set
of order $2,\ldots,b$. By the preceding estimates,
\[
 \Pp(\mathcal B^c)=O_{K,T}\bigl((\log n)^{-1}\bigr)
\]
and \eqref{eq:base-deletion-uncovered} holds uniformly on $\mathcal B$.
It remains to estimate $\mu(B_h)-\mu(B)$. 

If
$\mu(U_S\cap U_T)\ge b/(2C_0n)$, then \eqref{eq:uniform-bridge}
already follows from $b\asymp(\log n)^8$. Otherwise, $|S|>b$
and \eqref{eq:base-cap-mass} give
\begin{equation}\label{eq:exclusive-region-mass}
 \frac{b}{2C_0n}
 \le \mu(U_S)-\mu(U_S\cap U_T)
 =\mu(B)\le\frac12.
\end{equation}
Since $h\sqrt{d-1}=r/(\log n)^{3/2}\le1$, we may apply
Lemma~\ref{lem:spherical-growth}. Using \eqref{eq:radius-scale},
we obtain
\begin{equation}\label{eq:structure-scales}
\begin{aligned}
 \mu(B_h)-\mu(B)
 &\ge \frac{cr}{(\log n)^{3/2}}
        \left(\frac bn\right)^{1-1/(d-1)}\\
 &\ge \frac{c(\log n)^{1/(d-1)}b^{1-1/(d-1)}}
            {n(\log n)^{3/2}}\\
 &\ge \frac{c\sqrt b}{n(\log n)^{3/2}}
 \ge \frac{c(\log n)^{5/2}}{n},
\end{aligned}
\end{equation}
where $d\ge3$ and $b\asymp(\log n)^8$.

Substituting \eqref{eq:structure-scales} and
\eqref{eq:base-deletion-uncovered} into \eqref{eq:bridge-growth}
gives
\[
 \mu(U_S\cap U_T)
 \ge \frac{c(\log n)^{5/2}-C_{K,T}(\log n)^2}{n}
 \ge \frac{c'(\log n)^{5/2}}{n}
\]
for all sufficiently large $n$. The same event $\mathcal B$
works for every allowed deletion and every partition.
\end{proof}

\subsection{Independent reserve points merge the large components}
\label{sec:structure-merging}

Let the points $X_{M+1},\ldots,X_n$ form the reserve sample,
independent of the base.
For a fixed deletion set $Z\subseteq[n]$ with $|Z|\le K$,
list the surviving reserve points in increasing label order as
\[
 R_1,\ldots,R_N,
 \qquad N\ge n-M-K\ge\frac{n}{2\log n}.
\]
Conditional on the base set \(A\), these points remain independent and
uniform.

Take the components of $G[A\setminus Z]$ of order greater than
$b$ as the \emph{initial blocks}.
Let $Q_j$ be the number of connected groups of initial blocks
after exposing $R_1,\ldots,R_j$, using only edges from these
points to the initial blocks.
Thus $Q_j=1$ means that all initial blocks have merged.

Let $\mathcal E$ be the event that, for some $Z\subseteq[n]$ with
$|Z|\le K$, two initial blocks lie in distinct components of $G-Z$.
We show that this event is unlikely by proving that each
exposure contracts the conditional expectation of $Q_j-1$.

\begin{lemma}[Simultaneous merging of the large base components]
\label{lem:incidence}
For the base event $\mathcal B$ from Lemma~\ref{lem:bridge},
and all sufficiently large $n$,
\begin{equation}\label{eq:structure-merging-error}
 \Pp(\mathcal E\mid\mathcal B)
 \le C_Kn^{K+1}e^{-c(\log n)^{3/2}}
 =O_{K,T}\bigl((\log n)^{-1}\bigr).
\end{equation}
\end{lemma}

\begin{proof}

For the fixed deletion set $Z$, let
\[
 \mathcal G_j=\sigma(X_1,\ldots,X_M,R_1,\ldots,R_j),
 \qquad 0\le j\le N.
\]
Given $\mathcal G_j$, the current groups are determined and, for
$j<N$, the next reserve point $R_{j+1}$ is uniform on the sphere.
Work on the base event $\mathcal B$ and suppose $Q_j\ge2$.
A current group participates in a merger if the next point
meets both that group and another group.
Every current group is a union of initial blocks, so
Lemma~\ref{lem:bridge} gives conditional probability at least
\[
 \gamma=\frac{c(\log n)^{5/2}}n
\]
for each group to participate.
Summing these probabilities, the conditional expected number
of participating groups is at least $\gamma Q_j$.

All participating groups merge into one.
Since a merger involves at least two groups, the number of
groups lost is at least half the number participating.
Consequently,
\[
 \E(Q_j-Q_{j+1}\mid\mathcal G_j)
 \ge\frac{\gamma Q_j}{2},
 \qquad Q_j\ge2.
\]
When $Q_j=1$, all initial blocks have already merged and
$Q_{j+1}=1$. Hence, in every case,
\[
 \E(Q_{j+1}-1\mid\mathcal G_j)
 \le(1-\gamma/2)(Q_j-1).
\]
Iteration and Markov's inequality now give
\[
\begin{aligned}
 \Pp(Q_N>1\mid\mathcal G_0)
 &\le\E(Q_N-1\mid\mathcal G_0)\\
 &\le(Q_0-1)(1-\gamma/2)^N\\
 &\le n e^{-c(\log n)^{3/2}}.
\end{aligned}
\]

A union bound over at most $C_Kn^K$ deletion sets, followed by
averaging over base configurations conditional on $\mathcal B$,
proves \eqref{eq:structure-merging-error}.
\end{proof}

\begin{figure}[tbp]
\centering
\begin{tikzpicture}[
  font=\small, line cap=round, line join=round,
  group/.style={
    circle, draw=blue!60!black, fill=blue!12,
    minimum size=9mm, inner sep=1pt
  },
  bridge/.style={draw=orange!85!black, thick}
]
  \begin{scope}
    \node[group] at (0,.65)  {$S_1$};
    \node[group] at (2,.65)  {$S_2$};
    \node[group] at (1,-.75) {$S_3$};
    \node at (1,-1.75) {$Q_j$ groups};
  \end{scope}

  \draw[->,thick,black!65] (2.7,0)--(4.3,0)
    node[midway,above] {expose $R_{j+1}$};

  \begin{scope}[shift={(5.4,0)}]
    \filldraw[fill=blue!4,draw=blue!55!black]
      (1,0) ellipse[x radius=1.85,y radius=1.4];

    \node[group] (s1) at (0,.65)  {$S_1$};
    \node[group] (s2) at (2,.65)  {$S_2$};
    \node[group] (s3) at (1,-.75) {$S_3$};
    \coordinate (z) at (1,.08);

    \foreach \g in {s1,s2,s3}
      \draw[bridge] (z)--(\g);
    \fill[orange!85!black] (z) circle (2pt)
      node[below right=2pt,black] {$R_{j+1}$};

    \node at (1,-1.75) {$Q_{j+1}=Q_j-2$};
  \end{scope}
\end{tikzpicture}
\caption{One reserve point merges three current groups}
\label{fig:structure-incidence}
\end{figure}
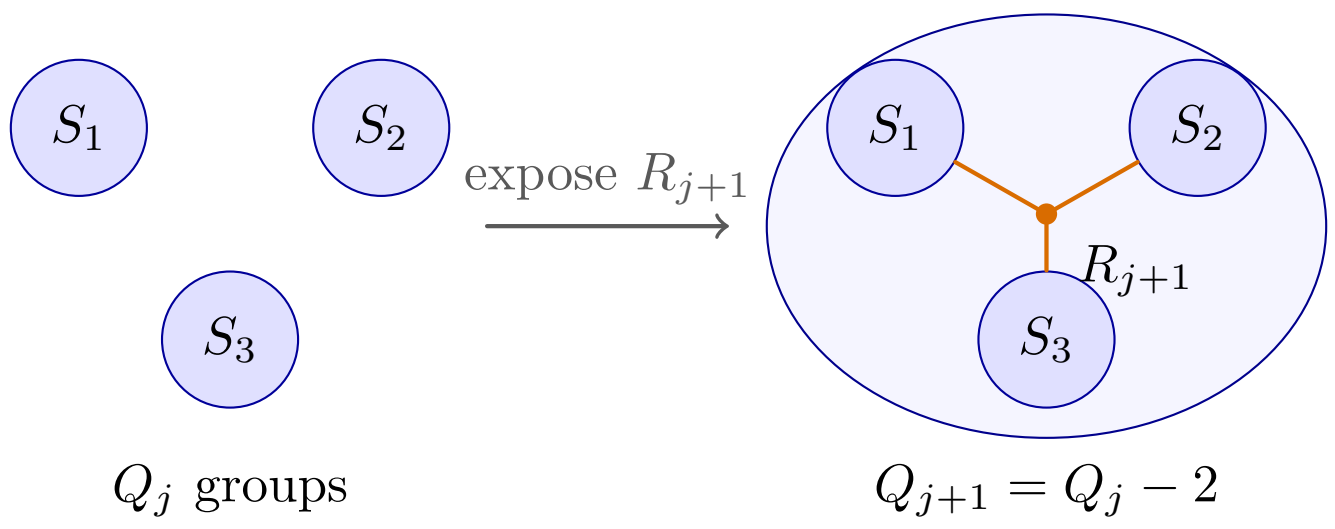

\subsection{The remaining vertices and completion of the proof}
\label{sec:structure-remainder}

On $\mathcal B\setminus\mathcal E$, for each deletion set
$Z\subseteq[n]$ with $|Z|\le K$, the component of $G-Z$
containing all components of $G[A\setminus Z]$ of order greater
than $b$ is called the \emph{principal component}.

\begin{proof}[Completion of the proof of Theorem~\ref{thm:robust}]
We exclude the following five bad events.
All probability bounds below are unconditional.

\begin{enumerate}
\item \emph{The base fails its required properties.}
Let $\mathcal F_1=\mathcal B^c$.
Lemma~\ref{lem:bridge} gives
\[
 \Pp(\mathcal F_1)
 =O_{K,T}\bigl((\log n)^{-1}\bigr).
\]

\item \emph{Merging fails despite a good base.}
Let $\mathcal F_2=\mathcal B\cap\mathcal E$.
By \eqref{eq:structure-merging-error},
\[
 \Pp(\mathcal F_2)=\Pp(\mathcal B)\Pp(\mathcal E\mid\mathcal B)=O_{K,T}\bigl((\log n)^{-1}\bigr).
\]

\item \emph{A reserve vertex has too few base neighbors.}
Let $\mathcal F_3$ be the event that some reserve vertex
has at most $K$ neighbors in the original base.
For each reserve label $M<i\le n$,
\[
 \deg_A(X_i)\sim\operatorname{Bin}(M,p).
\]
Indeed, conditional on $X_i$, the base points are independent,
and each is adjacent to $X_i$ with probability $p$.
The binomial lower-tail bound and \eqref{eq:base-window} give
\[
 \Pp\{\deg_A(X_i)\le K\}
 \le C_K(\log n)^Ke^{-Mp}
 \le\frac{C_{K,T}}n.
\]
A union bound over the $n-M=O(n/\log n)$ reserve labels yields
\begin{equation}\label{eq:reserve-base-neighbors}
 \Pp(\mathcal F_3)
 \le\frac{C_{K,T}(n-M)}n
 =O_{K,T}\bigl((\log n)^{-1}\bigr).
\end{equation}

\item \emph{Too many reserve vertices meet the exceptional base set.}
Let
\[
 R_Y=\#\{M<i\le n:X_i\in U_Y\},
 \qquad
 \mathcal F_4=\{R_Y>b/4\}.
\]
Thus $R_Y$ counts the reserve vertices adjacent to at least
one vertex of $Y$ in the original graph.
Conditional on the base, each reserve point lies in $U_Y$
with probability $\mu(U_Y)\le |Y|p$.
Consequently, \eqref{eq:base-low-degree} gives
\[
 \E R_Y
 =(n-M)\E\mu(U_Y)
 \le(n-M)p\,\E|Y|
 =O_{K,T}(1),
\]
since $(n-M)p=O(1)$.
Markov's inequality therefore gives
\begin{equation}\label{eq:structure-residual-error}
 \Pp(\mathcal F_4)
 \le\frac{4\E R_Y}{b}
 =O_{K,T}((\log n)^{-1}).
\end{equation}

\item \emph{The full graph contains a small separated set.}
Let $\mathcal F_5$ be the event that $G$ contains a
$K$-separated set of order between $2$ and $b$.
Since $b\le\lceil(\log n)^8\rceil$,
Proposition~\ref{prop:local}, applied with sample size $n$, gives
\[
 \Pp(\mathcal F_5)
 =O_{K,T}\bigl((\log n)^{-1}\bigr).
\]
\end{enumerate}

Suppose none of these events occurs, and fix any
$Z\subseteq[n]$ with $|Z|\le K$.
The principal component exists because $\mathcal B$ occurs
and $\mathcal E$ does not.
By the properties of $\mathcal B$, every surviving base vertex
outside this component belongs to $Y$, and $|Y|\le\log n$.

Since $\mathcal F_3$ does not occur, every surviving reserve
vertex retains at least one base neighbor after the deletion.
If it lies outside the principal component, that surviving
base neighbor must also lie outside, and hence belongs to $Y$.
Such a reserve vertex is therefore counted by $R_Y$.
Since $\mathcal F_4$ does not occur,
\[
\begin{aligned}
 \#\{\text{vertices outside the principal component of }G-Z\}
 &\le |Y|+R_Y\\
 &\le\log n+b/4<b.
\end{aligned}
\]

If a nontrivial component $S$ remained outside the principal
component, then
\[
 2\le|S|<b,\qquad
 G[S]\text{ is connected},\qquad
 \partial_G S\subseteq Z.
\]
Thus $S$ would be a $K$-separated set, contradicting the absence
of $\mathcal F_5$.
Every vertex outside the principal component is therefore
isolated. The principal component itself contains a base block
of order greater than $b$, so it is nontrivial.

Since $Z$ was arbitrary, this conclusion holds simultaneously
for every allowed deletion. Hence
\[
 \Pp\{\mathcal R_{K+1}(p)^c\}
 \le\sum_{i=1}^5\Pp(\mathcal F_i)
 =O_{K,T}\bigl((\log n)^{-1}\bigr),
\]
uniformly over $d\ge3$.
\end{proof}
\section{Minimum degree, connectivity, and hitting times}
\label{sec:consequences}

The robust structure theorem reduces connectivity to a
minimum-degree condition. We first estimate its probability,
then complete the fixed-radius connectivity and hitting-time
arguments.

\subsection{The minimum-degree probability}
\label{sec:low-degrees}

Let $X_a$ count vertices of degree exactly $a$. In the critical window
of Theorem~\ref{thm:main}, the mean in \eqref{eq:Lambda} satisfies
\begin{equation}\label{eq:low-degree-mean}
 \lambda_{K,n}=\E X_K
 =n\binom{n-1}{K}p^K(1-p)^{n-1-K}=\frac{e^{-t}}{K!}+o(1).
\end{equation}

\begin{theorem}[Minimum-degree probability]
\label{thm:low-degrees}
Fix an integer $K\ge0$ and $0\le T<\infty$.
Uniformly over $d\ge3$ and
\[
 np=\log n+K\log\log n+t,\qquad |t|\le T,
\]
we have
\begin{equation}\label{eq:minimum-degree-law}
 \Pp\{\delta(G_d(n,p))\ge K+1\}
 =e^{-\lambda_{K,n}}
   +O_{K,T}\bigl((\log n)^{-1}\bigr).
\end{equation}

\end{theorem}

\begin{proof}
Use the notation of Section~\ref{sec:geometry}.
The binomial degree distribution gives
\[
 \E X_a=O_{K,T}\bigl((\log n)^{a-K}\bigr),
 \qquad 0\le a<K.
\]
Consequently,
\[
 \left|\Pp\{\delta\ge K+1\}-\Pp(X_K=0)\right|
 \le\sum_{a=0}^{K-1}\E X_a
 =O_{K,T}\bigl((\log n)^{-1}\bigr),
\]
where the sum is empty when $K=0$.
It therefore suffices to prove
\begin{equation}\label{eq:low-degree-zero}
 \left|\Pp(X_K=0)-e^{-\lambda_{K,n}}\right|
 =O_{K,T}\bigl((\log n)^{-1}\bigr).
\end{equation}

\begin{enumerate}
\item \emph{Poissonization and a zero-count identity.}

First replace the sample size by an independent
$N\sim\Pois(n)$, keeping the radius $r$ fixed.
Let $\eta$ be the resulting point set and let $W$ count its
degree-$K$ vertices.
Poisson samples in disjoint regions are independent.

Fix $x\in\mathbb S^{d-1}$.
Independently sample $K$ uniform points $Y_1,\ldots,Y_K$ in $C_x$,
replace the points of $\eta$ inside $C_x$ by these points,
and insert $x$:
\[
 \eta^x=(\eta\cap C_x^c)\cup\{Y_1,\ldots,Y_K\}\cup\{x\}.
\]
Almost surely, $x$ has degree exactly $K$ in the graph on $\eta^x$:
its neighbors are precisely $Y_1,\ldots,Y_K$.
Figure~\ref{fig:local-resampling} illustrates this construction.

\begin{figure}[tbp]
\centering
\begin{tikzpicture}[
  font=\small,
  line cap=round,
  line join=round
]
  \foreach \shift/\name in {0/{\eta},6.6/{\eta^x}} {
    \begin{scope}[shift={(\shift,0)}]
      \filldraw[
        fill=blue!5,
        draw=blue!55!black,
        thick
      ] (0,0) circle (1.4);

      \node[blue!55!black] at (-.9,1.25) {$C_x$};
      \node at (0,2.15) {$\name$};

      \foreach \point in {
        (-1.8,.65),
        (-1.7,-1.2),
        (.7,1.75),
        (1.85,.2),
        (1.55,-1.5)
      } {
        \fill[blue!65!black] \point circle (2pt);
      }
    \end{scope}
  }

  \foreach \point in {
    (-.75,.45),
    (.15,.9),
    (.85,.25),
    (.5,-.65),
    (-.6,-.75)
  } {
    \fill[black!55] \point circle (2pt);
  }

  \draw[->,thick] (2.2,0)--(4.4,0);
  \node[align=center] at (3.3,.65)
    {Resample $K$ points\\inside $C_x$};
  \node at (3.3,-.4) {and insert $x$};

  \begin{scope}[shift={(6.6,0)}]
    \coordinate (x)  at (0,0);
    \coordinate (y1) at (.1,.9);
    \coordinate (y2) at (-.8,-.45);
    \coordinate (y3) at (.8,-.5);

    \foreach \point in {y1,y2,y3} {
      \draw[orange!85!black,thick] (x)--(\point);
      \fill[orange!85!black] (\point) circle (2.3pt);
    }

    \fill (x) circle (2.5pt);
    \node[above left] at (x) {$x$};
    \node[above] at (y1) {$Y_1$};
    \node[left]  at (y2) {$Y_2$};
    \node[right] at (y3) {$Y_3$};

    \node at (0,-2.05) {$\deg_{\eta^x}(x)=K$};
  \end{scope}
\end{tikzpicture}
\caption{Local resampling, illustrated for $K=3$.}
\label{fig:local-resampling}
\end{figure}
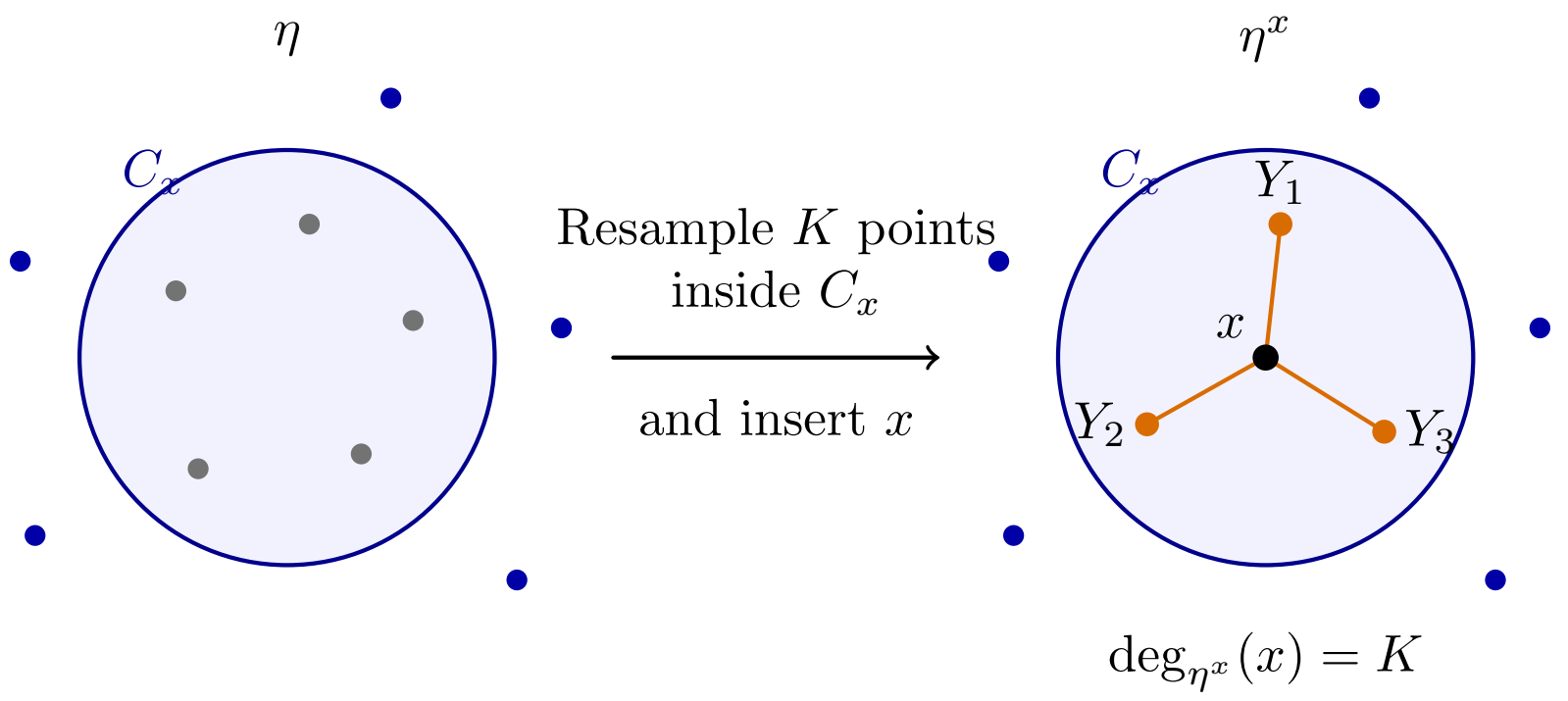

Let $W^x$ count the degree-$K$ vertices of $\eta^x$, including $x$;
in particular, $W^x\ge1$.

Set
\[
 g(s)=\E s^W=\Pp(W=0)+\sum_{j=1}^{\infty}\Pp(W=j)s^j.
\]
Distinguishing one point of the Poisson sample, using
$j\Pp(N=j)=n\Pp(N=j-1)$ and rotational invariance, gives
\[
 \lambda_*:=\E W=\frac{ne^{-np}(np)^K}{K!},
 \qquad
 g'(s)=\lambda_*\E s^{W^x-1},
 \qquad 0<s<1.
\]
For nonnegative integers $a,b$ and $0\le s\le1$,
$|s^a-s^b|\le|a-b|$.
Since $g(0)=\Pp(W=0)$ and $g(1)=1$, integration yields
\begin{equation}\label{eq:zero-coupling}
\begin{aligned}
 \left|\Pp(W=0)-e^{-\lambda_*}\right|
 &=\left|\int_0^1e^{-\lambda_*s}
          \bigl(g'(s)-\lambda_*g(s)\bigr)\,ds\right|\\
 &\le(1-e^{-\lambda_*})\E|W-(W^x-1)|\\
 &\le\E|W-(W^x-1)|.
\end{aligned}
\end{equation}

\item \emph{The effect of resampling one cap.}

After subtracting the forced contribution of $x$, the discrepancy
$|W-(W^x-1)|$ is bounded by the sum of three counts:
removed degree-$K$ points, added degree-$K$ points other than $x$,
and original points outside $C_x$ whose degree-$K$ indicators
change.
We estimate their expectations separately.

First, by rotational invariance, the expected number of
degree-$K$ points in a region is its measure times $\lambda_*$.
Since all original points in $C_x$ are removed and $\mu(C_x)=p$,
the expected number of removed degree-$K$ points is $p\lambda_*$.

Next, consider the added points other than $x$.
For $K=0$, there are none, so assume $K\ge1$.
For an added point $y\in C_x$, the forced point $x$ is a neighbor,
and all original points in $C_y\setminus C_x$ remain.
Their number has distribution $\Pois(n[p-q(x,y)])$.
Thus $y$ can have degree $K$ only if
\[
 \Pois\bigl(n[p-q(x,y)]\bigr)\le K-1.
\]
Put $u=\dist(x,y)/r\le1$. By \eqref{eq:cap-union},
\[
 p-q(x,y)=\mu(C_x\cup C_y)-p\ge pu/\pi.
\]
Also, \eqref{eq:chord-density} gives
$F'(s)=a_d s\bigl(s\sqrt{1-s^2/4}\bigr)^{d-3}$.
Since $s\sqrt{1-s^2/4}$ is nondecreasing on $[0,r]$ for
$r\le\sqrt2$, $F$ is convex there. Hence
$F(ur)\le uF(r)=up$.
For any $a\ge0$,
\[
 \Pp\{\Pois(a)\le K-1\}
 =e^{-a}\sum_{j=0}^{K-1}\frac{a^j}{j!}
 \le C_K e^{-a/2}.
\]
Consequently,
\[
 \Pp\{\Pois(n[p-q(x,y)])\le K-1\}
 \le C_K\exp\left(-\frac{npu}{2\pi}\right).
\]
Each of the $K$ added points is uniform in $C_x$, with
distribution $\mu(dy)/p$.
Their expected degree-$K$ count is therefore at most
\[
\begin{aligned}
 &\frac Kp\int_{C_x}
   \Pp\{\Pois(n[p-q(x,y)])\le K-1\}\,\mu(dy)\\
 &\qquad\le
   \frac{C_K}{p}\int_0^1
       e^{-npu/(2\pi)}\,dF(ur)\\
 &\qquad\le
   C_K\left(
       e^{-np/(2\pi)}
       +\frac{np}{2\pi}\int_0^1u e^{-npu/(2\pi)}\,du
       \right)\\
 &\qquad=
   \frac{2\pi C_K}{np}\left(1-e^{-np/(2\pi)}\right)
 \le\frac{C_K}{np}.
\end{aligned}
\]
The second inequality follows by integration by parts, using
$F(0)=0$, $F(r)=p$, and $F(ur)\le pu$.

Finally, for any $K\ge0$, consider an original point
$y\notin C_x$ and write $q=q(x,y)$.
Its degree-$K$ indicator can change only if both of the following
occur: there are at most $K$ other original points in
$C_y\setminus C_x$, and at least one removed or added point
lies in $C_y\cap C_x$.
These two events are independent.

Since $\dist(x,y)>r$, \eqref{eq:cap-union} gives
\[
 p-q=\mu(C_x\cup C_y)-p\ge p/\pi.
\]
The number of other original points in $C_y\setminus C_x$
has distribution $\Pois(n[p-q])$, so the first event has
probability at most
\[
 C_K(np)^K e^{-np/\pi}.
\]
The expected numbers of removed and added points in
$C_y\cap C_x$ are $nq$ and $Kq/p$, respectively.
Thus the second event has probability at most
\[
 (n+K/p)q.
\]
Summing over the original points outside $C_x$, and using
$\int q(x,y)\,\mu(dy)=p^2$ from Lemma~\ref{lem:overlap},
bounds the expected number of changed
indicators by
\[
\begin{aligned}
 &C_Kn(n+K/p)(np)^K e^{-np/\pi}
       \int_{C_x^c}q(x,y)\,\mu(dy)\\
 &\qquad\le C_K(np)^{K+2}e^{-np/\pi}.
\end{aligned}
\]

Combining the three contributions gives
\[
 \E|W-(W^x-1)|
 \le p\lambda_*+\frac{C_K}{np}
       +C_K(np)^{K+2}e^{-np/\pi}
 =O_{K,T}\bigl((\log n)^{-1}\bigr).
\]
Together with \eqref{eq:zero-coupling}, this proves the required
zero-count estimate for the Poisson sample.

\item \emph{Return to a fixed sample size.}

Realize $\eta$ as the first $N$ points of an infinite uniform
sample independent of $N$.
Let $X_K^{(j)}$ be the degree-$K$ count on the first $j$ points,
always using the original radius $r$.
Then $X_K^{(N)}$ has the law of $W$, while $X_K^{(n)}=X_K$.

Adding one point changes the count only through the new vertex
and adjacent old vertices of degree $K$ or $K-1$.
Consequently,
\[
\begin{aligned}
 \E|X_K^{(j+1)}-X_K^{(j)}|
 &\le\Pp\{\operatorname{Bin}(j,p)=K\}\\
 &\quad+jp\left[
   \Pp\{\operatorname{Bin}(j-1,p)=K\}
   +\Pp\{\operatorname{Bin}(j-1,p)=K-1\}\right].
\end{aligned}
\]
The term involving $K-1$ is omitted when $K=0$.
Put $m=\lceil\sqrt n\log n\rceil$.
Uniformly for $|j-n|\le m$, we have $jp=np+o(1)$, so the
preceding bound gives
\[
 \E\bigl|X_K^{(j+1)}-X_K^{(j)}\bigr|
 \le C_{K,T}\frac{\log n}{n}.
\]

We split according to whether $N$ lies within $m$ of $n$:
\[
 \Pp\{X_K^{(N)}\ne X_K^{(n)}\}
 \le \Pp\{|N-n|>m\}+\Pp\{|N-n|\le m,\ X_K^{(N)}\ne X_K^{(n)}\}.
\]
Since $N\sim\Pois(n)$, a Poisson tail bound gives
\[
 \Pp\{|N-n|>m\}
 \le 2e^{-c(\log n)^2}.
\]
For the second term, a single-step estimate gives
\[
 \Pp\{|N-n|\le m,\ X_K^{(N)}\ne X_K^{(n)}\}
 \le \sum_{j=n-m}^{n+m-1}
       \E\bigl|X_K^{(j+1)}-X_K^{(j)}\bigr|
 \le C_{K,T}\frac{(\log n)^2}{\sqrt n}.
\]
Combining the two bounds yields
\[
 \Pp\{X_K^{(N)}\ne X_K^{(n)}\}
 =o\bigl((\log n)^{-1}\bigr).
\]
In particular,
\[
 \bigl|\Pp(W=0)-\Pp(X_K=0)\bigr|
 =o\bigl((\log n)^{-1}\bigr).
\]

It remains to compare the two means. By \eqref{eq:low-degree-mean}
and the formula for $\lambda_*$,
\[
 \frac{\lambda_{K,n}}{\lambda_*}
 =\left[\prod_{j=1}^{K}\left(1-\frac jn\right)\right]
   e^{np}(1-p)^{n-1-K}
 =1+O_{K,T}\left(\frac{(\log n)^2}{n}\right),
\]
where the last step uses $\log(1-p)=-p+O(p^2)$ and
$p=O(\log n/n)$. Since $\lambda_{K,n}=O_{K,T}(1)$ by
\eqref{eq:low-degree-mean}, it follows that
\[
 \left|e^{-\lambda_*}-e^{-\lambda_{K,n}}\right|
 \le |\lambda_*-\lambda_{K,n}|
 =O_{K,T}\left(\frac{(\log n)^2}{n}\right)
 =o\bigl((\log n)^{-1}\bigr).
\]
Combining this with the preceding estimates gives
\[
\begin{aligned}
 \left|\Pp(X_K=0)-e^{-\lambda_{K,n}}\right|
 &\le \left|\Pp(X_K=0)-\Pp(W=0)\right|
      +\left|\Pp(W=0)-e^{-\lambda_*}\right|\\
 &\quad+\left|e^{-\lambda_*}-e^{-\lambda_{K,n}}\right|
 =O_{K,T}\bigl((\log n)^{-1}\bigr),
\end{aligned}
\]
which proves \eqref{eq:low-degree-zero}.
\end{enumerate}
\end{proof}

\subsection{Proofs of the main theorems}

\begin{proof}[Proof of Theorem~\ref{thm:main}]
Write $G=G_d(n,p)$. On $\mathcal R_{K+1}(p)$, deleting at most
$K$ vertices leaves one nontrivial component and possibly isolated
vertices. If $\delta(G)\ge K+1$, every surviving vertex retains
a neighbor, so the remaining graph is connected. Together with
the inclusions
\[
 \{G\text{ is $(K+1)$-vertex-connected}\}
 \subseteq
 \{G\text{ is $(K+1)$-edge-connected}\}
 \subseteq
 \{\delta(G)\ge K+1\},
\]
this shows that all three events coincide on $\mathcal R_{K+1}(p)$.

By Theorem~\ref{thm:robust}, the probability that these events
differ is $O_{K,T}(1/\log n)$. Theorem~\ref{thm:low-degrees}
therefore gives
\[
 e^{-\lambda_{K,n}}+O_{K,T}(1/\log n)
\]
for both connectivity probabilities. The limiting laws follow from
$\lambda_{K,n}=e^{-t}/K!+o(1)$. Taking $K=0$ in
Theorem~\ref{thm:robust} gives \eqref{eq:ordinary-structure}.
\end{proof}

The $O_{K,T}(1/\log n)$ error uses the exact mean $\lambda_{K,n}$.
For $K\ge1$, replacing it by $e^{-t}/K!$ can introduce an error
of order $(\log\log n)/\log n$.

\begin{proof}[Proof of Theorem~\ref{thm:hitting}]
Fix $T>0$ and write
\[
 p_\pm=\frac{\log n+K\log\log n\pm T}{n},
 \qquad G_\pm=G_d(n,p_\pm).
\]
 For each pair $i<j$, define
\[
 A_{ij}=\{ij\in E(G_+)\setminus E(G_-)\},
 \qquad
 B_{ij}=\{|N_{G_-}(i)\cup N_{G_-}(j)|\le K\}.
\]
Thus $A_{ij}$ says that the edge $ij$ appears during the window,
and $B_{ij}$ says that at its left endpoint the two vertices
have at most $K$ distinct neighbors between them.

We claim that
\[
 \{p_{\delta,K}\ne p_{\kappa,K}\}
 \cap\{p_{\delta,K}\in[p_-,p_+]\}
 \cap\mathcal R_{K+1}(p_-)
 \subseteq
 \bigcup_{i<j}(A_{ij}\cap B_{ij}).
\]
Indeed, on the event on the left, the graph at $p_{\delta,K}$ has minimum
degree $K+1$, yet some deletion set $Z$ of at most $K$ vertices
disconnects it. At $p_-$, the graph $G_--Z$ had one nontrivial
component and otherwise only isolated vertices. That component
remains connected at $p_{\delta,K}$. Every surviving vertex has a neighbor
at $p_{\delta,K}$, so a component disjoint from the original one must
contain an edge $ij$. Both $i$ and $j$ were isolated in $G_--Z$.
Hence $ij$ is a new edge, giving $A_{ij}$, and
\[
 N_{G_-}(i)\cup N_{G_-}(j)\subseteq Z,
\]
giving $B_{ij}$. This proves the inclusion.

To estimate the right-hand side, first note that
$\Pp(A_{ij})=p_+-p_-$. On $A_{ij}$, the endpoints are farther
apart than $r_d(p_-)$. Applying \eqref{eq:cap-union} at this radius,
the union of their two caps has mass at least $(1+1/\pi)p_-$.
Conditional on the endpoint locations, the other $n-2$ points
are independent and uniform, and those in the union are precisely
the neighbors counted in $B_{ij}$. Thus
\[
 \Pp(B_{ij}\mid A_{ij})
 \le \Pp\{\operatorname{Bin}(n-2,(1+1/\pi)p_-)\le K\}.
\]
A union bound now gives
\[
\begin{aligned}
 \Pp\!\left(\bigcup_{i<j}(A_{ij}\cap B_{ij})\right)
 &\le \sum_{i<j}\Pp(A_{ij})\Pp(B_{ij}\mid A_{ij})\\
 &\le \binom n2(p_+-p_-)
      \Pp\{\operatorname{Bin}(n-2,(1+1/\pi)p_-)\le K\}\\
 &\le C_{K,T}\,n(np_-)^K e^{-(1+1/\pi)np_-}
 =O_{K,T}(n^{-1/\pi}).
\end{aligned}
\]

By the claimed inclusion and Theorem~\ref{thm:robust},
\[
\begin{aligned}
 \Pp\{p_{\delta,K}\ne p_{\kappa,K}\}
 &\le \Pp\{p_{\delta,K}\notin[p_-,p_+]\}
      +\Pp\{\mathcal R_{K+1}(p_-)^c\}
      +\Pp\!\left(\bigcup_{i<j}(A_{ij}\cap B_{ij})\right)\\
 &\le \Pp\{\delta(G_-)\ge K+1\}
      +\Pp\{\delta(G_+)<K+1\}+o(1),
\end{aligned}
\]
uniformly over $d\ge3$. Applying Theorem~\ref{thm:low-degrees}
at the two endpoints yields
\[
 \limsup_{n\to\infty}\sup_{d\ge3}
 \Pp\{p_{\delta,K}\ne p_{\kappa,K}\}
 \le
 \exp\!\left(-\frac{e^T}{K!}\right)
 +1-\exp\!\left(-\frac{e^{-T}}{K!}\right).
\]
Letting $T\to\infty$ shows that the two hitting times coincide
with probability tending to one. Since
$p_{\delta,K}\le p_{\mathrm e,K}\le p_{\kappa,K}$,
the same holds for all three.
\end{proof}
\appendix
\section{Cap unions, intersections, and spherical neighborhood growth}
\label{app:cap-geometry}
\subsection{Union and intersection of two caps}\label{app:two-cap-proof}

\begin{proof}[Proof of Lemma~\ref{lem:crescent}]
We first give the union bound \eqref{eq:cap-union}.
The claim is immediate if $x=y$. If $y=-x$, the two caps have
disjoint interiors and their union has measure $2p$.
Otherwise let $P=\operatorname{span}\{x,y\}$ and let $\Pi$ be the
orthogonal projection onto $P$. For a uniform point $Z$ on the
sphere, write
\[
 \Pi Z=\rho W,\qquad \rho=\|\Pi Z\|.
\]
Almost surely $\rho>0$, and rotational symmetry makes $W$ uniform
on the unit circle in $P$ conditional on $\rho$. Since $x,y\in P$,
\[
 Z\in C_x\iff x\cdot W\ge\tau/\rho,\qquad
 Z\in C_y\iff y\cdot W\ge\tau/\rho.
\]
For $\rho<\tau$ both events are empty. Otherwise they correspond
to arcs $A_x,A_y$ of the same half-length
$\alpha=\arccos(\tau/\rho)\le\pi/2$, as in
Figure~\ref{fig:crescent-slices}.

Let $\theta$ be the shorter arc distance from $x$ to $y$.
The two outward half-arcs together contribute $2\alpha$.
Between the centers, the arcs cover a further length
$\min\{\theta,2\alpha\}$. Thus, on each slice with $\alpha>0$,
\[
 \frac{|A_x\cup A_y|}{|A_x|}
 =1+(\frac{\theta}{2\alpha}\wedge1)
 \ge 1+\frac{u\wedge1}{\pi}.
\]
Indeed, $\theta\ge\dist(x,y)=ur$, and
$\alpha\le\arccos\tau=2\arcsin(r/2)\le\pi r/2$.
Averaging the arc lengths over $\rho$ proves \eqref{eq:cap-union}.

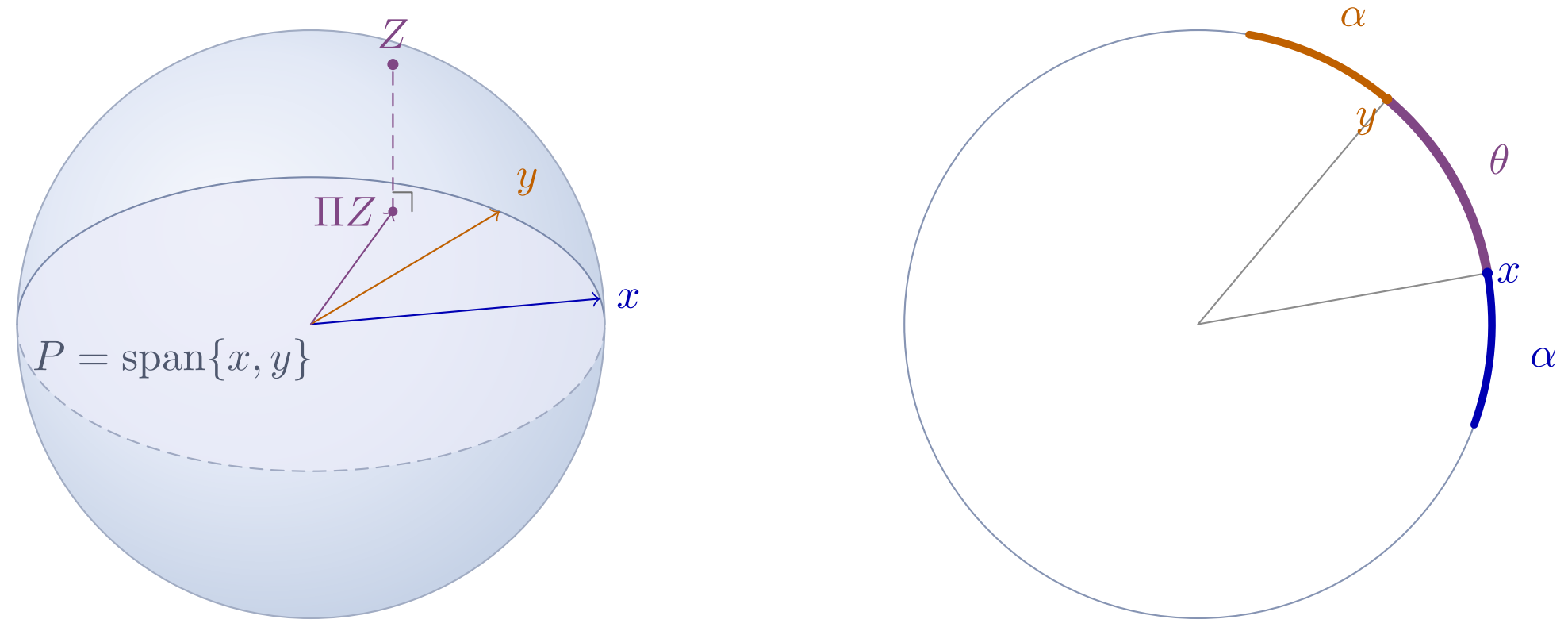
\begin{figure}[tb]
\centering
\begin{tikzpicture}[scale=2.45,font=\small,line cap=round,line join=round]
 \definecolor{arcGlobeLight}{HTML}{F5F7FC}
 \definecolor{arcGlobeMid}{HTML}{E3E9F6}
 \definecolor{arcGlobeShade}{HTML}{BECCE2}
 \definecolor{arcPlane}{HTML}{EBECF9}
 \definecolor{arcPlaneEdge}{HTML}{7988AA}
 \definecolor{arcGain}{HTML}{804785}
 \pgfdeclareradialshading{crescentSphere}{\pgfpoint{-8bp}{9bp}}{
   color(0bp)=(arcGlobeLight);
   color(17bp)=(arcGlobeMid);
   color(33bp)=(arcGlobeShade);
   color(50bp)=(arcGlobeShade)}
 \tikzset{
   x mark/.style={blue!70!black},
   y mark/.style={orange!75!black},
   new mark/.style={arcGain}
 }

 \begin{scope}
   \pgfmathsetmacro{\rhoPlot}{sqrt(2/3)}
   \pgfmathsetmacro{\projX}{\rhoPlot*cos(70)}
   \pgfmathsetmacro{\projY}{.5*\rhoPlot*sin(70)}
   \coordinate (X) at ({cos(10)},{.5*sin(10)});
   \coordinate (Y) at ({cos(50)},{.5*sin(50)});
   \coordinate (Q) at (\projX,\projY);
   \coordinate (Z) at (\projX,{\projY+.5});

   \shade[shading=crescentSphere] (0,0) circle (1);
   \fill[arcPlane,fill opacity=.69] (0,0)
     ellipse[x radius=1,y radius=.5];
   \draw[arcPlaneEdge!72!white,densely dashed] (-1,0)
     arc[start angle=180,end angle=360,x radius=1,y radius=.5];
   \draw[arcPlaneEdge] (1,0)
     arc[start angle=0,end angle=180,x radius=1,y radius=.5];
   \draw[arcPlaneEdge!70!white] (0,0) circle (1);

   \draw[x mark,->] (0,0)--(X);
   \draw[y mark,->] (0,0)--(Y);
   \draw[new mark,->] (0,0)--(Q);
   \draw[new mark,densely dashed] (Q)--(Z);
   \draw[black!55] ({\projX+.065},\projY)
     -- ({\projX+.065},{\projY+.065})
     -- (\projX,{\projY+.065});
   \fill[new mark] (Q) circle (.016) node[left] {$\Pi Z$};
   \fill[new mark] (Z) circle (.019) node[above] {$Z$};
   \node[x mark,right] at (X) {$x$};
   \node[y mark,above right] at (Y) {$y$};
   \node[arcPlaneEdge!65!black] at (-.47,-.12) {$P=\operatorname{span}\{x,y\}$};
   \node[below] at (0,-1.09) {(a) Projection onto \(P\)};
 \end{scope}

 \begin{scope}[shift={(3.02,0)}]
   \draw[arcPlaneEdge!90!white] (0,0) circle (1);
   \draw[black!45,thin] (0,0)--(10:1) (0,0)--(50:1);

   \draw[x mark,line width=1.8pt] (-20:1)
     arc[start angle=-20,end angle=10,radius=1];
   \draw[new mark,line width=2.2pt] (10:1)
     arc[start angle=10,end angle=50,radius=1];
   \draw[y mark,line width=1.8pt] (50:1)
     arc[start angle=50,end angle=80,radius=1];

   \fill[blue!70!black] (10:1) circle (.018) node[right,inner sep=2pt] {$x$};
   \fill[orange!75!black] (50:1) circle (.018) node[below left,inner sep=2pt] {$y$};
   \node[x mark,right] at (-6:1.08) {$\alpha$};
   \node[new mark,right] at (31:1.09) {$\theta$};
   \node[y mark,above right] at (66:1.05) {$\alpha$};

   \node[below] at (0,-1.09) {(b) Splitting the arc union at \(x,y\)};
 \end{scope}
\end{tikzpicture}
\caption{From spherical caps to circular arcs.}
\label{fig:crescent-slices}
\end{figure}

Then we give the overlap bound \eqref{eq:far-overlap}.
To estimate the threshold $\tau$, observe that the density of
$x\cdot Y$ is bounded by $a_d$ because $d\ge3$. By symmetry,
\[
 \frac12-p=\Pp\{0\le x\cdot Y<\tau\}\le a_d\tau.
\]
The assumption $p\le 2\log n/n$ and $a_d\asymp\sqrt d$ then gives $\tau\sqrt d\ge c$
for an absolute constant $c>0$ and all sufficiently large $n$.

For $0\le u<1$, recall that
\[
 J(u)=\int_u^1 (1-z^2)^{(d-3)/2}\,dz.
\]
The substitution $t=((1-z^2)/(1-u^2))^{(d-1)/2}$ gives
\begin{equation}\label{eq:Ju}
 J(u)=\frac{(1-u^2)^{(d-1)/2}}{d-1}
 \int_0^1\frac{dt}{\sqrt{1-(1-u^2)t^{2/(d-1)}}}.
\end{equation}
Note that
\[
 1-(1-u^2)t^{2/(d-1)}\ge 1-(1-u^2)=u^2
\]
and
\[
\begin{aligned}
 \int_0^1 \frac{dt}{\sqrt{1-(1-u^2)t^{2/(d-1)}}}
 &\ge \left(\int_0^1
       \bigl[1-(1-u^2)t^{2/(d-1)}\bigr]\,dt\right)^{-1/2}\\
 &=\frac{1}{\sqrt{u^2+2(1-u^2)/(d+1)}}.
\end{aligned}
\]
 We obtain
\begin{equation}\label{eq:mills}
 \frac{(1-u^2)^{(d-1)/2}}
 {(d-1)\sqrt{u^2+2(1-u^2)/(d+1)}}
 \le J(u)\le
 \frac{(1-u^2)^{(d-1)/2}}{(d-1)u},
 \qquad 0<u<1.
\end{equation}
At $u=\tau$, the two denominators are comparable because
$\tau\sqrt d\ge c$. Using $p=a_dJ(\tau)$ and $a_d\asymp\sqrt d$
gives
\[
 p\asymp\frac{(1-\tau^2)^{(d-1)/2}}{\tau\sqrt d}.
\]

Taking logarithms gives
\[
 \log(1/p)
 =\frac{d-1}{2}\log\frac1{1-\tau^2}
  +\frac12\log(d\tau^2)+O(1)
 \asymp d\log\frac1{1-\tau^2},
\]
hence
\begin{equation}\label{eq:d-tau^2}
d\tau^2
\asymp d\wedge\left(d\log\frac1{1-\tau^2}\right)
\asymp d\wedge\log(1/p)
\asymp d\wedge\log n.
\end{equation}

 If $y=-x$, then $C_x\cap C_y=\varnothing$
 and \eqref{eq:far-overlap} holds. Otherwise, denote $s=\dist(x,y)$ and $u=2\tau/\|x+y\|
=\tau/\sqrt{1-s^2/4}$. On the unit sphere,
\[
\begin{aligned}
 C_x\cap C_y
 &=\{z:x\cdot z\ge\tau,\ y\cdot z\ge\tau\}\\
 &\subseteq \{z:(x+y)\cdot z\ge2\tau\}\\
 &=\left\{z:\frac{x+y}{\|x+y\|}\cdot z\ge u\right\}.
\end{aligned}
\]
If $u\ge1$, the last set has measure zero and \eqref{eq:far-overlap}
again holds. For $u<1$, the last set is a cap of measure $a_dJ(u)$
by rotation invariance. Thus we obtain
\[
 \frac{q(x,y)}{p}
 =\frac{\mu(C_x\cap C_y)}{\mu(C_x)}
 \le \frac{a_dJ(u)}{a_dJ(\tau)}
 =\frac{J(u)}{J(\tau)}.
\]

Since $u\ge\tau$, the integrand in \eqref{eq:Ju} at threshold $u$
is at most the integrand at threshold $\tau$ for every $t\in[0,1]$.
Thus the ratio of the integral factors is at most one, and
\[
 \frac{J(u)}{J(\tau)}
 \le \left(\frac{1-u^2}{1-\tau^2}\right)^{(d-1)/2}
 = \left(1-\frac{u^2-\tau^2}{1-\tau^2}\right)^{(d-1)/2}
 \le \exp\left(-\frac{d-1}{2}
                   \frac{u^2-\tau^2}{1-\tau^2}\right).
\]
Note that
\[
 u^2-\tau^2=\frac{\tau^2s^2}{4-s^2}\ge\frac{\tau^2s^2}{4},
 \qquad
 1-\tau^2=r^2(1-r^2/4)\le r^2.
\]
Substituting these bounds gives
\[
 \frac{q(x,y)}p
 \le \exp\left(-\frac{(d-1)\tau^2}{8}(s/r)^2\right).
\]
When $s\ge r/10$, we have $(s/r)^2\ge1/100$. Using \eqref{eq:d-tau^2}, it follows that
\[
 \frac{q(x,y)}{p}
 \le \exp\!\left(-\frac{(d-1)\tau^2}{800}\right)
 \le e^{-a(d\wedge\log n)}
\]
for an absolute constant $a>0$. This proves \eqref{eq:far-overlap}.
\end{proof}

\subsection{Spherical neighborhood growth}\label{app:growth-proof}

\begin{proof}[Proof of Lemma~\ref{lem:spherical-growth}]
By spherical isoperimetry
\cite[Theorem~5.1.6]{Vershynin2018}, among sets of the same
measure, spherical caps minimize the measure of the chordal
$h$-neighborhood.
Let $C$ be a cap with $\mu(C)=\mu(B)$ and angular radius
$0<\theta\le\pi/2$.
Comparing the integral factors in \eqref{eq:Ju} at thresholds
$\cos\theta$ and $0$ gives
\[
\begin{aligned}
\mu(C)
&=a_d\frac{(\sin\theta)^{d-1}}{d-1}
  \int_0^1\frac{dt}{\sqrt{1-\sin^2\theta\,t^{2/(d-1)}}}\\
&\le a_d\frac{(\sin\theta)^{d-1}}{d-1}
  \int_0^1\frac{dt}{\sqrt{1-t^{2/(d-1)}}}\\
&=\frac12(\sin\theta)^{d-1}
\end{aligned}
\]

 For $0\le t\le h$, we have
\[
 \sin(\theta+t)\ge\sin\theta\cos h,
 \qquad
 (\cos h)^{d-2}
 \ge 1-\frac{(d-2)h^2}{2}\ge\frac12.
\]
Since chordal distance does not exceed geodesic distance,
$C_h$ contains the cap of angular radius $\theta+h$.
Thus the added spherical band yields
\[
\begin{aligned}
 \mu(B_h)-\mu(B)
 &\ge \mu(C_h \setminus C)\\
 &\ge a_d\int_0^h\sin^{d-2}(\theta+t)\,dt\\
 &\ge \frac{a_dh}{2}(\sin\theta)^{d-2}\\
 &\ge c h\sqrt{d-1}\,\mu(B)^{1-1/(d-1)},
\end{aligned}
\]
where we used $a_d\asymp\sqrt{d-1}$ in the last step.
\end{proof}

\section*{AI assistance}

In 2025, the author proved the connectivity result in the
regime $d=\Omega(\log^3 n)$. Working from the author's original
manuscript, ChatGPT 6.0 suggested ways to strengthen the argument.
With its assistance, the author completed the proof presented here.
The author takes full responsibility for the mathematical
arguments and the final text.

\end{document}